\documentclass[10pt]{amsart}
\usepackage{amsmath}
\usepackage{amsxtra}
\usepackage{amscd}
\usepackage{amsthm}
\usepackage{amsfonts}
\usepackage{amssymb}
\usepackage{mathrsfs}
\usepackage{color} 

\newtheorem{lem}{Lemma}[section]
\newtheorem{cor}[lem]{Corollary}

\newtheorem{prop}[lem]{Proposition}

\newtheorem{ax}[lem]{Axiom}

\newtheorem{defn}[lem]{Definition}

\newtheorem{rem}[lem]{Remark}

\theoremstyle{definition}

\theoremstyle{remark}

\newcommand{\thmref}[1]{Theorem~\ref{#1}}
\newcommand{\secref}[1]{\S\ref{#1}}
\newcommand{\lemref}[1]{{ Lemma~\ref{#1}}}
\newcommand{\defref}[1]{Definition~\ref{#1}}
\newcommand{\propref}[1]{{ Proposition~\ref{#1}}}
\newcommand{\corref}[1]{Corollary~\ref{#1}}

\newcommand{\remref}[1]{Remark~\ref{#1}}

\def\eqn#1#2{ \begin{equation} \label{#1} #2 \end{equation} }
\def\eq#1{\begin{equation*}  #1 \end{equation*}}

\newcommand{\nc}{\newcommand}
\nc{\renc}{\renewcommand}
\nc{\ssec}{\subsection}
\nc{\sssec}{\subsubsection}
\nc{\on}{\operatorname}
\nc{\remm}[1]{\<{remark} \ \lbl{#1} \>{remark}}

\renc\k{\mathbf{k}}
\def\th#1{{\theta_{#1}}}
\nc\ol{\overline} 
\nc\wh{\widehat}
\nc\tboxtimes{\wt{\boxtimes}}

\nc{\Aa}{{\mathbb{A}}}
 \nc{\Gg}{{\mathbb{G}}}

\def\g{\gamma}

\def\a{\alpha}  \def\b{\beta} \def\i{\iota}
\def\Vol{{\rm Vol}}

\nc{\Hh}{{\mathbb{H}}}

\nc{\Nn}{{\mathbb{N}}}
\nc{\Pp}{{\mathbb{P}}}
\nc{\Rr}{{\mathbb{R}}}
\nc{\BV}{{\mathbb{V}}}
\nc{\BW}{{\mathbb{W}}}
\nc{\Zz}{{\mathbb{Z}}}
\nc{\Qq}{{\mathbb{Q}}}
\nc{\Ss}{{\mathbb{S}}}
\nc{\Cc}{{\mathbb{C}}}
\def\FF{\mathcal F}
 
\def\om{\omega}
\nc{\Oo}{{\mathcal{O}}}
\nc{\Mm}{{\mathcal{M}}}

\nc{\eo}{{\mathbf{\tau}}}
\nc{\dU}{{\overset{\bullet}{\bigcup}}{}}
\nc{\du}{\, {\overset{.}{\cup}}\, {}}
\nc{\dual}[1]{{\overset{\vee}{#1}{}}}

\nc{\cM}{{\check{\mathcal M}}{}}
 \nc{\oM}{{\overset{\circ}{\mathcal M}}{}}

 \nc{\fB}{{\mathfrak{B}}}

\nc{\tT}{{\widetilde{T}}}  

\nc{\fF}{{\mathcal{F}}}

\nc{\bb}{{\mathbf{b}}}
\nc{\bc}{{\mathbf{c}}}
\nc{\bd}{{\mathbf{d}}}
\nc{\be}{{\mathbf{e}}}
\nc{\bj}{{\mathbf{j}}}
\nc{\bn}{{\mathbf{n}}}

\nc{\bph}{{\mathbf{\phi}}}
\nc{\bp}{{\mathbf{p}}}
\nc{\bq}{{\mathbf{q}}}
\nc{\bF}{{\mathbf{F}}}
\nc{\bu}{{\mathbf{u}}}
\nc{\bv}{{\mathbf{v}}}
\nc{\bx}{{\mathbf{x}}}
\nc{\bh}{{\mathbf{h}}}
\nc{\bs}{{\backslash}}
\nc{\by}{{\mathbf{y}}}
\nc{\bw}{{\mathbf{w}}}
\nc{\bA}{{\mathbf{A}}}
\nc{\bK}{{\mathbf{K}}}
\nc{\bI}{{\mathbf{I}}}
\nc{\bB}{{\mathbf{B}}}
\nc{\bG}{{\mathbf{G}}}
\nc{\bC}{{\mathbf{C}}}
\nc{\bD}{{\mathbf{D}}}
\nc{\bP}{{\mathbf{P}}}
\nc{\bH}{{\mathbf{H}}}
\nc{\bM}{{\mathbf{M}}}
\nc{\bN}{{\mathbf{N}}}
\nc{\bV}{{\mathbf{V}}}
\nc{\bU}{{\mathbf{U}}}
\nc{\bL}{{\mathbf{L}}}
\nc{\bT}{{\mathbf{T}}}
\nc{\bW}{{\mathbf{W}}}
\nc{\bX}{{\mathbf{X}}}
\nc{\bY}{{\mathbf{Y}}}
\nc{\bZ}{{\mathbf{Z}}}
\nc{\bS}{{\mathbf{S}}}
\nc{\ba}{{\mathbf{a}}}

\nc{\sA}{{\mathsf{A}}}
\nc{\sB}{{\mathsf{B}}}
\nc{\sC}{{\mathsf{C}}} 
\nc{\sF}{{\mathsf{F}}}
\nc{\sG}{{\mathsf{G}}}
\nc{\sK}{{\mathsf{K}}}
\nc{\sM}{{\mathsf{M}}}
\nc{\sO}{{\mathsf{O}}}
\nc{\sQ}{{\mathsf{Q}}}
\nc{\sP}{{\mathsf{P}}}
\nc{\sZ}{{\mathsf{Z}}}
\nc{\sfp}{{\mathsf{p}}}
\nc{\sr}{{\mathsf{r}}}
\nc{\sg}{{\mathsf{g}}}
\nc{\sff}{{\mathsf{f}}}
\nc{\sfb}{{\mathsf{b}}}
\nc{\sfc}{{\mathsf{c}}} 

\nc{\tA}{{\widetilde{{A}}}}
 \def\RV{{\rm RV}}

\def\m{\setminus}
\nc{\tD}{{\widetilde{{A}}}}
\nc{\tH}{{\widetilde{{A}}}}
\nc{\tB}{{\widetilde{{B}}}}
 
\nc{\tG}{{\widetilde{G}}}
\nc{\TM}{{\widetilde{\mathbb{M}}}{}}
\nc{\tO}{{\widetilde{\mathsf{O}}}{}} 
\nc{\TZ}{{\tilde{Z}}}
\nc{\tx}{{\tilde{x}}}
\nc{\tf}{{\tilde{f}}}
\nc{\tz}{{\tilde{\zeta}}}
\nc{\tmu}{{\tilde{\mu}}}
\nc{\td}{{\tilde{d}}}
        
\nc{\tX}{{\widetilde{X}}}
  
   \nc{\E}{{\mathop{\operatorname{\rm E }}}}
 \nc{\Mor}{{\mathop{\operatorname{\rm Mor \,}}}}
\nc{\Ob}{{\mathop{\operatorname{\rm Ob \,}}}}
  \nc{\Sym}{{\mathop{\operatorname{\rm Sym}}}}
   \nc{\Aut}{{\mathop{\operatorname{\rm Aut}}}}
 \nc{\Spec}{{\mathop{\operatorname{\rm Spec}}}}
\nc{\Ker}{{\mathop{\operatorname{\rm Ker}}}}
 \nc{\dom}{{\mathop{\operatorname{\rm dom}}}}
\nc{\End}{{\mathop{\operatorname{\rm End}}}}
 \nc{\Hom}{\on{\Hom}} 
 \nc{\GL}{{\mathop{\operatorname{\rm GL}}}}
 \nc{\Id}{{\mathop{\operatorname{\rm Id}}}}
 \nc{\rk}{{\mathop{\operatorname{\rm rk}}}}
\nc{\irk}{{\mathop{\operatorname{\rm i-rk}}}}
 \nc{\length}{{\mathop{\operatorname{\rm length}}}}
\nc{\supp}{{\mathop{\operatorname{\rm supp}}}}
\nc{\val}{{\rm val}}

\nc{\valr}{{\rm val_{rv}}}
\nc{\valrv}{\valr}
\def\bdd{^{\rm bdd}}
\def\Gb{{\G_A \bdd}}
\def\Gbt{{\G_{A(t)} \bdd}}

 \nc{\Res}{{\rm res}}
 \nc{\Aff}{{\rm Aff}}
\nc{\res}{{\mathop{\operatorname{\rm res}}}}
\nc{\rad}{{\mathop{\operatorname{\rm rad}}}}

\def\tensor{{\otimes}}

\def\meet{\cap}
\def\union{\cup}

\def\G{\Gamma}
\def\<{\begin}
 \def\>{\end}

\nc{\tV}{{\widetilde{{V}}}}
\nc{\hb}[1]{\hbox{#1}}
\def\rv{{\rm rv}}

\def\Pi{\prod}
 
\nc{\seq}[1]{\stackrel{#1}{\sim}}

\nc{\oeq}[1]{\underset{#1}{=}} 
\def\eg{\oeq{gr}}
\def\inv {{^{-1}}}
 \def\beq#1{\begin{equation} \label{#1}}   
\def\eeq{\end{equation}}

\def\iso{\simeq}
\def\RR{{\mathcal R}}
\def\prf{\begin{proof}}
\def\eprf{\end{proof} }

 \def\lbl#1{  \label{#1}  }

\def\RR{{\mathcal R}}

\usepackage{color}

\nc{\Claim}[1]{{\noindent \bf Claim{ #1 }}}
  \nc{\pr}{{\mathop{\operatorname{\rm pr}}}}

\nc{\Mmm}{{(1+\Mm)}}
  
   \def\RVp{{\RV^{>0}}}

\def\RES{{\bf RES}}

  \nc{\SG}{{\mathop{\operatorname{{  K_+ }}}}} 
    \nc{\SGo}{{\mathop{\operatorname{{  K_{+0} }}}}} 
   \nc{\K}{{\mathop{\operatorname{ { K} }}}}

\def\col#1#2{\begin{pmatrix} #1 \\ #2 \end{pmatrix} }

\def\ax{{\a x}}

\newcommand{\ga}{\gamma}
\newcommand{\de}{\delta}
\newcommand{\ep}{\epsilon}

\newcommand{\Ga}{\Gamma}

\newcommand{\Om}{\Omega}

\newcommand{\cA}{\mathcal A}

\newcommand{\cH}{\mathcal H}

\newcommand{\cO}{\mathcal O}

\newcommand{\bR}{\mathbb R}

\newcommand{\sI}{\mathscr I}

\newcommand{\hM}{\widehat{M}}

\newcommand{\QQ}{\qquad\quad}

\newcommand{\bmat}{\left ( \begin{matrix} }
\newcommand{\emat}{\end{matrix} \right ) }

\nc{\q}{{\mathfrak{q}}}
\nc{\qq}{\dot{q}}

\nc{\Fnc}{{\rm Fn_{*}}}
\nc{\Fncb}{{\rm Fn_{*}^{bdd}}}
 
\def\volG{{{\rm vol } \G}}
\def\vol{{\rm vol}}
 \def\Var{{\rm Var}}
 \def\VF{{\rm VF}}
  \def\V{{\rm V}}

 \def\RG{{\mathcal R}^{\G}}
\def\RGb{{\mathcal R}^{\G \bdd}}
\def\RGi{{\mathcal R}^{\G \infty}}
 
 \author{Ehud Hrushovski, David Kazhdan}
 
\address{\newline Institute of Mathematics, the Hebrew
 University of Jerusalem, Givat Ram, Jerusalem, 91904, Israel.} 
 
 \email{ ehud@math.huji.ac.il,kazhdan@math.huji.ac.il}

\title{The value ring of geometric motivic integration, and the Iwahori Hecke algebra of $SL_2$}

\begin{document}

\maketitle

\vspace{-0.5cm}
\centerline{(with an appendix by Nir Avni)}
\vspace{0.5cm}

 \<{section}{Introduction}

 In \cite{HK}, an integration theory for valued fields was developed with a Grothendieck group approach.   Two types of categories were studied.  The first was of semi-algebraic sets over a valued field, with all semi-algebraic morphisms.  The Grothendieck  ring of this category
 was shown to admit two natural homomorphisms, esssentially into the Grothendieck  ring of varieties over the residue field.  These can be viewed as generalized Euler characteristics.
The objects of the second category are semi-algebraic sets with  volume forms; the 
morphisms are semi-algebraic bijections preserving the absolute value of the volume form.
(Some finer variants were also studied.)  The Grothendieck ring of bounded objects in this category can be viewed as a universal integration theory.  

Even before the restriction to bounded sets, an   isomorphism was shown between the semiring  of semi-algebraic sets with  measure preserving morphisms, and certain semirings formed out twisted varieties over
the residue field, and rational polytopes over the value group.   Though this description
is very precise, the target remains complicated.  
With a view to representation-
theoretic applications, we require a simpler description of the  possible values of the integration,
and in particular natural homomorphisms into fields.  In the present paper we obtain such results after tensoring
with $\Qq$,  in particular 
 introducing additive inverses.  Since this operation trivializes the full semiring, we restrict to bounded sets.   We show that the resulting $\Qq$-algebra is generated by its one-dimensional part.  
In the ``geometric'' case, i.e.
working over an elementary submodel as a base, we determine  the structure 
precisely.  As a corollary we obtain useful
canonical homomorphisms in the general case.  

Let $F$ be a valued field of residue characteristic $0$.  Let $V$ be an $F$-variety.
A {\em semialgebraic subset } of $V$ is a Boolean combination of subvarieties and of
sets defined by valuation inequalities $\{x \in U: \val f(x) \leq \val g(x) \}$, where 
$U$ is a relatively closed $F$-subvariety of $V$, and 
$f,g$ are regular functions on $U$.   (It is possible to think of the $F^a$-points 
defined by these equalities, but better to think of $K$-points where $K$ is an undetermined
valued field extension of $F$.)    

 Let $\Vol_F$ be the category of   semi-algebraic sets with bounded semi-algebraic volume forms;
 see \ref{vol} for a precise definition.  The Jacobian of any semi-algebraic map between such objects  can then be defined, outside
 a lower dimensional variety; morphisms are  semi-algebraic bijections whose Jacobian
 has valuation zero (outside a lower dimensional variety.)  
The Grothendieck ring  $\K(\Vol_F)$ of this category can be viewed as a universal integration theory
for semialgebraic sets and volume forms over $F$.  This ring is graded by dimension,
but one can form out of it a ring $\K^{df}(\Vol_F)$ of ``pure numbers'', ratios of integrals of equal dimension (see \S 1.1).  We state there a version of \thmref{vf} in the case of a higher dimensional local field.

Let  $\Var_{\bF}$ be the category of algebraic varieties over the residue field of $F$. 
  $\K^{df}_\Qq(\Var_{\bF})$ is the dimension-free Grothendieck ring with rational
coefficients this category.  There exists a natural homomorphism $L_F^{Var} : \K^{df}_\Qq(\Var_{\bF}) \to 
 \K^{df}(\Vol_F)$, induced by taking the full pullback of a variety $V \subseteq \Aa^n(\bF)$ to the valuation ring, with the 
standard form $dx_1\ldots dx_n$.   
  
 Assume $F$ has value group generated by $n$ elements $\g_1,\ldots,\g_n$. 
  Extend $L_F^{Var}$  to a  homomorphism 
$$ L_F: \K^{df}_{\Qq} (\Var_{\bF}) [t_1,\ldots,t_n,q_1,\ldots,q_n] \to  \K^{df}(\Vol_F)$$
by mapping $q_i$ to the ratio of the annulus  of valuative radius $\g_i$ to the unit annulus $U_0$; and  $t_i$ to the logarithmic quantity
  $L_F(t_i) = [(\{x: 0 \leq \val(x) < \g_i \}), dx / x) ] / [(U_0, dx)]$. 
 
Localizations by certain elements will be needed.  They are explained in the text
 before the statement of \thmref{vf}.    Here we will just denote them with a subscript $loc$.
We denote by $L_F$ the homomorphism induced on localizations also.

\<{thm} \lbl{vf1}  Assume   $F$ has value group $\Zz^n$.  Let $\bF$ denote the residue field of
  $F$.  

There exists a canonical homomorphism
 $$I_F: 
 \K^{df}(\Vol_F)_{loc} \to  \K^{df}_{\Qq} (\Var_{\bF}) [t_1,\ldots,t_n,q_1,\ldots,q_n]_{loc}$$
 with $I_F J_F = Id$.  
  \>{thm}

  The ring $\K^{df}_{\Qq}(\Var_{\bF})$ is a subring of the usual Grothendieck ring $K_\Qq(\Var_{\bf})$ of varieties over $F$, localized at $[G_m]$;
  we have $\K^{df}_{\Qq}(\Var_{\bF}) \cong \{a/[G_m]^k: k \in \Nn, a \in K_\Qq(\Var_{\bf}), \dim(a) \leq k \}$.  
it includes an element $L=1+\frac{[pt]}{[G_m]}$, corresponding on the left    
  to  the ratio  of the volume of a closed and an open ball of the same radius.   On the other hand, $q_i$ corresponds to the ratio
  of the ball $\val(x)>\val(t_i)$, to the unit ball $\val(x)>0$.    
 
   { The quantities $L,q_1,\ldots,q_n$ are
  $\Qq$-algebraically independent.  }   This contrasts with the $p$-adic integration theories,
  and   those of Denef, Denef-Loeser, Cluckers-Loeser, where one has ($n=1$ and)  
 $L \inv=q_1$.  The reason for the additional degree of freedom is that we chose the ``geometric'' realization
  of the universal integral.   It can already be seen via the  the following functoriality in ramified extensions:
 
   If $F \leq F'$ is a finite ramified field extension, whose value group is generated (for simplicity)
 by $\g_1/m_1,\ldots,\g_n/m_n$, then we have:
  $$I_{F'}: 
 \K^{df}(\Vol_{F'})_{loc} \to  \K^{df}_{\Qq} (\Var_{\bF'})_{loc} [t'_1,,\ldots,t'_n,q'_1,\ldots,q'_n]_{loc}$$
 
 With $m_i t'_i = t_i$ and $(q'_i)^{m_i} = q_i$.  At the limit over all ramified extensions,
 or just a family whose value groups approach $\Qq^n$,  the homomorphisms
 $I_{F'}$ become an {\em isomorphism}.    In fact the fundamental case here is really the case of  divisible  value group.
 
 Viewed as an integral, $I_F$ satisfies Fubini and the usual change of variable formula,
 with respect to arbitrary semi-algebraic maps.  It is also additive with respect
 to definable maps into the value group or residue field.    
 
  In the 
 case of value group $\Zz^n$ described above, the theorem should be compared to
 earlier integration theories of Fesenko and Parshin ; see  \cite{fesenko}.

The above statements are all special cases of the results in \cite{HK}, with 
improvement only in the description of the target ring.  This  
 depends  on a closer study of the Grothendieck ring of
bounded piecewise linear polytopes.  We express in closed form the motivic volume
of any bounded polytope over an ordered Abelian group,
in terms of quantities $\i(b)$ referring to the length of a one-dimensional segment
$[0,b)$, and Boolean quantities $e(b)$ that can be viewed as referring to the existence
or not of $b$ as a  rational point.  Note that  $\frac{1}{m}\i(x) \neq \i(\frac{x}{m})$ in general.  
 The formulas specialize (in their graded version) to standard integration formulas, and
on the other hand formulas giving the number of integer points in bounded polytopes.  
But since they must also be valid in groups such as $\Zz^n$, nothing can be assumed
about the index of arithmetic sequences.  Nevertheless when sufficient care is taken
with arithmetic issues, it turns out that the formulas can be 
    proved using
integration by parts.

In \cite{HK}, a parallel theory without volume forms, and without ignoring lower dimensional sets, was also developed.  On the one hand, a universal invariant was found, with 
values in a Grothendieck ring formed out of $K(\Var_{\bF})$ and $K(\G)$.  
(Theorem 1.1)
On the
other hand, two homomorphisms were found, essentially into $K(\Var_{\bF})$;
they were deduced from the universal invariant and  two ``Euler characteristic'' homomorphisms $K(\G) \to \Zz$, found earlier
by \cite{marikova} and \cite{kageyama-fujita}. (Theorem 10.5)   However, no universality property
was shown for the latter.  
The two Euler characteristics are known to be universal with respect to
$GL_n(\Qq)$ transformations, but it is 
 $GL_n(\Zz)$ transformations that are relevant here; since it is these (along with
 translations by values of rational points) that  lift to 
the valued field.         \thmref{unbdd}  fills this gap   in the rational coefficient case, by showing
that even with respect to integral transformations alone, $\K^{df}(\G) \cong \Qq^2$. 

In the appendix we define the Iwahori Hecke algebra of $SL_2$ over an algebraically closed valued field. Iwahori Hecke algebras are usually defined for (quasi-)split algebraic groups over non archimedian local fields as convolution algebras with respect to the Haar measure. Here, instead, we use motivic integration. We give an analogue of the Bernstein presentation for the algebra and find its center. In \cite{Lee}, a construction of the Iwahori Hecke algebra of $SL_2$ over a two dimensional local field is given. We think this construction is unrelated to ours.
 
{\bf Acknowledgment} The authors were partially supported by ISF grants \# 244/03 and 1461/05.
 
\>{section}
\<{section}{The Grothendieck ring of bounded polytopes over an ordered Abelian group}
\ssec{The dimension-free part of a graded ring}

While we are ultimately interested in $\Qq$-algebras,  in the interest of simpler proofs
we will also use semirings for the basic lemmas.   Elements of the Grothendieck semiring are represented by definable sets, and equality corresponds to definable bijections.  For the corresponding ring representing an element $[X]-[Y]$ requires two definable sets, and equality $[X]-[Y] = [X']-[Y']$ invokes a third
definable set $Z$ and an isomorphism $X \du Y' \du Z \to X' \du Y \du Z$.  Thus 
a canonical  isomorphism between semirings, when available, is not only stronger but easier to prove
than the isomorphism of rings it implies. 

All semigroups in this paper will be commutative, with addition denoted by $+$, and a distinguished
element $0$ (perhaps the term {\em commutative monoid} is more standard.)

Given a graded   semiring $R = \oplus_{n \geq 0} R_n$, and an  element $a_1 \in R_1$,
 $R[a_1 \inv]$ is naturally $\Zz$-graded; let 
  $ R^{df}_{a_1}  = R[{a_1}  \inv]_0$ be the zero'th homogeneous component.  When $a_1$
  is fixed we will just write $R^{df}$.
  We think of the elements of $R^{df}$ as ratios or pure numbers, whereas the elements
  of $R$ may have ``units''.

 As a semigroup,  $R^{df}_a$ can also be described as the direct limit of the semigroups $R_d$ under the maps $R_d \to R_{d+1}$ given by $x \mapsto a_1 x$.   In some cases that
 will be encountered, 
 e.g. when $R_d$ is the Grothendieck group of varieties of dimension $\leq d$, $R^{df}$
 can be thought of as a stabilized version of the Grothendieck group of varieties (of all dimensions at once.)
 
Define a semiring homomorphism $f:R \to  R[a \inv] _0$   by $f(r) = \frac{r}{{a_1}^n}$
for $r \in R_n$.  $R[a_1 \inv]_0$ has the universal property for semiring homomorphisms
$g: R \to S$ such that $g(a_1)=1$. 

The  Laurent polynomial semiring $R^{df} [t,t \inv]$ is
isomorphic, as a
 $\Zz$- graded semiring, to the localization $R[a _1\inv]$ (with $t \mapsto a_1$.)
 
 If $f_i: A \to B_i$ is a semiring homomorphism, $B_1 \tensor_A B_2$ is defined to be the
universal semiring $B$ with maps $g_i: B_i \to B$ such that $g_1f_1 = g_2f_2$.  If
$\bA,\bB,\bB_i$ are the ring canonically obtained from $A,B,B_i$ by introducing additive inverses,
one verifies immediately that the natural map $\bB \to (\bB_1 \tensor_A \bB_2)$ is an isomorphism. 
 
 \<{lem} \lbl{gr}  Let $\phi: R_1 \tensor R_2 \to R_3 $ be a surjective
  homomorphism of graded
semirings   with $\phi(e_1 \tensor 1) = \phi(1 \tensor e_2) =  e_3$, with  kernel $\sim $.  Let $S_i = R_i[e_i \inv]_0$.

If  $\sim $   is generated by the single relation
$1 \tensor e_2 \sim e_1 \tensor 1$, then $\phi$ induces
an  isomorphism $S_1 \tensor S_2 \to S_3$. 

More generally, if $\sim $  is generated by $1 \tensor e_2 = e_1 \tensor 1$ and  $1 \tensor f_2 = f_1 \tensor 1$, with $e_i,f_i \in R_i[1]$,
then  $\phi$ induces 
a surjective  homomorphism $S_1 \tensor S_2 \to S_3$, with kernel generated by $\frac{f_2}{e_2} \tensor 1 \sim 1 \tensor \frac{f_1}{e_1}$.
 \>{lem}

\prf  
We pass to the localizations, and obtain 
 a  homomorphism of Laurent polynomial semirings
\[  S_1[t_1, t_1 \inv] \tensor S_2 [t_2, t_2 \inv] )   \to S_3 [t_3,t_3 \inv] \] 
restricting to a homomorphism $S_1 \tensor S_2 \to S_3$,
with kernel generated by $ 1 \tensor t_2 \sim t_1 \tensor 1 $ and  (in the second case) an additional relation, that may be written
$t_2 \frac{f_2}{e_2} \tensor 1 \sim 1 \tensor \frac{f_1}{e_1} t_1$
 But 
$S_1[t_1, t_1 \inv] \tensor S_2 [t_2, t_2 \inv]  / (1 \tensor t_2 = t_1 \tensor 1) =
(S_1 \tensor S_2) (t, t\inv)$.     We thus have a surjective homomorphism 
$(S_1  \tensor S_2) [t , t  \inv] ) \to S_3[t,t \inv]$ with kernel generated by $   (\frac{f_2}{e_2} \tensor 1)t \sim (1 \tensor \frac{f_1}{e_1}) t $,
or equivalently by   $   \frac{f_2}{e_2} \tensor 1 \sim 1 \tensor \frac{f_1}{e_1} 1 $.
Restricting to $S_1 \tensor S_2$ we find a   homomorphism into $S_3$; it is easy to see that it must be surjective, with kernel generated by the same relation. \eprf

\<{lem} \lbl{gr2}  Let $R$ be a graded ring, $a_1 \in R_1$, $R^{df} = R^{df}_{a_1}$.  
Let  $b \in R_1$, $I=Rb$, $\bR = R/I$,
$\ba_1 = a_1 / I \in \bR$.  Let $I^{df} = R^{df} \frac{b}{a_1}$.  
Then 
$R^{df} / I^{df} \cong \bR ^{df}$.
\>{lem}

\prf The homomorphism $R \to R/I$ extends to a homomorphism 
$h: R[ a_1 \inv] \to \bR [ \ba_1 \inv]$ of $\Zz$-graded rings.  $h$
is surjective on every homogeneous component.  In particular $h$ restricts
to a surjective ring homomorphism $h_0: R[a_1 \inv]_0 \to \bR[\ba_1 \inv]_0$.
Any element of $R[a_1 \inv]_0$ can be written as $\frac{r}{a_1^n}$ for some
$r \in R_n$.  If $h_0(\frac{r}{a_1^n})  =0$ then $h(r) \ba_1^m=0$ for some $m \geq 0$.
So $h(ra_1^m)=0$, i.e. $ra_1^m = bs$ for some $s$.  Since $r,a_1,b$ are homogeneous
of repsective degrees $n,1,1$, we  can take $s$ to be homogeneous
of degree $n+m-1$.  But then in $R[a_1 \inv]_0$ we have 
$\frac{r}{a_1^n} = \frac{b}{a_1} \frac{s}{a_1^{n+m-1}} \in I^{df}$.  This shows
that $\ker (h_0) = I^{df} $, proving the lemma.  
\eprf

We also have:

\<{lem} \lbl{gr3}  Let $R,S$ be graded semirings, $e \in R_1, e' \in S_1$,  and let 
$f: R \to S$ be an injective homomorphism, $f(e)=e'$.  If for any $r' \in S$, for some $n$, 
$r' (e')^n \in f(R)$, then $f$ induces an isomorphism $R^{df}_{e} \to S^{df}_{e'}$. \>{lem}

\prf Clear.  \eprf

\ssec{Two categories of bounded definable subsets of $\G^n$}
 
    Throughout the text, $A$ denotes an ordered Abelian group, seen as a base
 subset of a model of the theory $DOAG$ of  divisible  ordered Abelian groups.

 \<{defn} \lbl{Gcat} 
 
 (1)  An object of  $ \G_A[n]$ is a   subset  of $\G^n$ defined by linear equalities and inequalities
with $\Zz$-coefficients and parameters in $A$.   When $A$ is fixed, we write $\G[n] = \G_A[n]$.
Given
$X,Y \in \Ob \G[n]$, $f \in \Mor_\G(X,Y)$ iff $f$ is a bijection, and there exists a partition $X=\union_{i=1}^n X_i$,
$M_i \in \GL_n(\Zz), a_i \in A^n$, such that for $x \in X_i$,
$$f(x) = M_i x + a_i$$  

 (2)   $\Gb [*]$ is the full subcategory of $\G [*]$
consisting of bounded sets, i.e. an element of  $\Ob \Gb [n]$ is a definable subset of $[-\g,\g]^n$
for some $\g \in \G$.

(3)    $\Ob {\volG}_A[n] = \Ob \G_A[n]$   Given
$X,Y \in \Ob {\volG}_A[n]$,  $f \in \Mor_{\volG _A[n]}(X,Y)$ iff $f \in \Mor_{\G[n]}$ and
for any $x=(x_1,\ldots,x_n) \in X$, if $y=(y_1,\ldots,y_n) = f(x)$ then $\sum_{i=1}^n x_i = \sum_{i=1}^n y_i$.  

 (4)   $\volG _A \bdd [n] = \Ob \Gb [n]$ is the 
 full subcategory of $\volG [n]$  with objects $X \subseteq [\g,\infty)^n$
 for some $\g \in \G$.     (Such objects will be called semi-bounded.)

(5) ${\volG_A}[*]$ is the direct sum of the categories ${{\volG}}[n]$ over $n \geq 0$; similarly for the other categories.

\>{defn}
 
$\SG[\Gb] [n]$ denotes the Grothendieck semigroup of $\Gb [n]$. By definition, it is the free
semigroup generated by the objects of $\Gb [n]$, subject to the relations:
$[X_1]+[X_2]=[Z]$ when there exists a partition $Z = Z_1 \du Z_2$ of $Z$, $Z_i \in \Gb [n]$,
with $X_i,Z_i$ isomorphic in $\Gb [n]$.  The zero object is defined to be the class of $\emptyset$.  
  It is easy to see (using boundedness) that for $n \geq 1$, 
$\SG \Gb[n]$ has finite direct sums (represented by disjoint unions).  Hence
any element of $\SG \Gb [n]$ is represented by an object of $\Gb [n]$.  

The semigroup $\SG \Gb[0]$ is $\Nn$; in this case only $0$ and $1$ are 
represented by an object of $\Gb [0]$.

 $\SG \Gb $
is the graded semiring $\oplus_{n \in \Nn} \SG \Gb[n]$.   Here $\SG \Gb[0] = \Nn$.   $\K \Gb$ is the corresponding
ring. Similar notation is used for the measured categories.  
 
Observe that a disjoint union of $\volG[n]$  isomorphisms is again a $\volG$
 isomorphism, provided that it is a $\G[n]$ isomorphism. 

 Here we will be interested in dimension-free quantities, i.e. ratios of elements
 of $\G[n]$ for each $n$, taking their direct limit over $n$.    
 We will normalize $\SG[\Gb]$ using the element $[0]_1$.  Let 
 $$\SGo(\Gb) = ( \SG[\Gb] [  [0]_1 \inv])_0$$
Let $\K^{df}(\Gb)$ be the corresponding ring, and 
 \
  $$ \K_{\Qq}^{df} (\Gb) = \Qq \tensor \K^{df}(\Gb)$$
 
\<{remark} \label{dehomogenizer}  We defined the dimension free ring using a dehomgenizing element $a=[0]_1$,
but could define a variant $\SG^{df(b)}$ using   $b=[X]_1$, for any nonempty definable $X \subseteq \G$.  The choice $a=[0]_1$ has the following universality property:  $\SG^{df(b)}(\Gb)$   embeds into a localization  
$\SG^{df} (\Gb)[ (\frac{[X]}{[0]_1}) \inv] $  
of $\SG^{df}(\Gb)$, $\frac{c}{b^n} \mapsto \frac{c}{a^n} (\frac{b}{a})^{-n}$.     This requires a lemma:
if $[Y \times \{0\}^m] = [Y' \times  \{0\}^m] $ then $[Y \times X^m ] = [Y' \times X^m]$.   \>{remark}  

For $a \in \Qq \tensor A$, let $e(a) = [a]_1 / [0]_1 $.  We have $[a]_1[0]_1 = [(a,0)]_2 = [(a,a)]_2 = [a]_1^2$,
using the $GL_2(\Zz)$ map $(x,y) \mapsto (x+y,y)$.  Hence $e(a)$ is   idempotent.
For $a \in A$, we have $e(a) = 1$.    We also have
 an element $\i(a) = [0,a)_1 / [0]_1 $, in $\K^{df}(\Gb)$.   (Here $[0,a)$ is the closed-open interval, for $a>0$; if $a<0$ we let $\i(a)=-\i(-a)$, and $\i(0)=o$.)    We will sometimes  write $[a,b)$ to denote the class $\i(b)-\i(a)$.
 If $\phi(x_1,\ldots,x_n)$ is a formula, we will sometimes write $[\phi]$ for
 the class of $\{(x_1,\ldots,x_n): \phi(x_1,\ldots,x_n) \}$.
 
 We define the dimension-free Grothendieck ring as in the unmeasured case:
$$\K_{+}^{df} (\vol \Gb) = ( \SG[{\vol \Gb} ][*][  [0]_1 \inv])_0$$

Let $\K^{df}(\vol \Gb)$ be the corresponding ring, and 
  $$ \K_{\Qq}^{df} (\vol \Gb) = \Qq \tensor \K^{df}(\vol \Gb)$$

It turns out that the measured ring can be constructed   from the unmeasured
one; we thus  begin by studying the latter. 

\ssec{Definable functions}  \lbl{definable-functions}
Recall the semigroup of functions $Fn(\G, \SG (\Gb))$.   An element of this semiring is represented
by a definable set $F \subseteq \G \times \G^m$, such that $F(x) = \{y: (a,y) \in F \}$
is bounded for any $x$.  $F$ represents a function in the following sense:  given
any ordered Abelian group extension $A(t)$ of $A$, generated over $A$ by a single element $t$, we obtain an element $[F(t)]$ of $\SG( \Gbt)$.

Similarly we define $Fn(\G, \SGo (\Gb))$.  An element is again represented by 
 a definable set $F \subseteq \G \times \G^m$, such that $F(x) = \{y: (a,y) \in F \}$
is bounded for any $x$.  
Two such sets $F,F'$ represent the same function if
for any $A'$ extending $A$ and $b \in A'$, 
$[F(b)] /[0]_1^m= [F'(b)] /[0]_1^m$ as elements of $\K^{df}(\G \bdd _{A'})$, i.e.
if $[F(b)]_{m+m'} = [{F'(b)}]_{m+m'}$ for some $m'$.     
Note that $e(t)$ represents the function $1$ in this formalism, since
$[b]_1=[0]_1$ in $\K^{df}(\G \bdd _{A'})$, using the translation $x \mapsto x-b$.   Hence $F(t)$,
 $e(t)F(t)$ represent the same function.  Since all $A'$ are at issue, we may take $A'=A(b)$.
  Addition is
 defined pointwise on representatives.  There is more than one option for multiplication;
 at present we will use pointwise 
multiplication, yielding a semiring.    
 The ring of functions $Fn(\G,\K^{df}(\Gbt) = \Zz \tensor Fn(\G, \SGo (\Gb))$ is the ring of formal differences;
an element $[F_1]-[F_2]$   is represented by a pair $(F_1,F_2)$, with the obvious
rules for equivalence, sum and product. 

If $F(x)$ represents an element of $Fn(\G, \SGo (\Gb))$, and $h: \G \to \G$ is any definable
function, consider $[F \circ h] \in Fn(\G, \SGo (\Gb))$.  If $[F]=[F']$ then 
$[e(h(x))] [F \circ h ] = [e(h(x))] [F' \circ h]$.  In particular, if  $h(x) = nx + a$, with $a \in A$
and $n \in \Nn$, then $[F]=[F']$ implies $[F \circ h] =  [F' \circ h]$.  But  if $h$ has non-integral coefficients, this need not be the case.   

\ssec{Integral notation}

Let $f$ be a function represented by $F$.
If $a<b \in \G$, 
write $\int_{a}^{b} f(x) dx $ for the class of $\{(t,y): a \leq t < b, (t,y) \in F \}$.
Note that $\int_{a}^b f(x) dx =  \int_{a}^b f(x) e(x) dx $.

One can think of the element ``dx'' as denoting the idempotent $e(x)$. 

If $a>b$, we let $\int_{a}^b dt = - \int_{b}^a dt$.

If $\a \in \Qq$ and $c \in \Qq \tensor A$,
we have a term  $e(\a t - c) \in Fn(\G,\K(\Gb))$, mapping $b$ to   the idempotent $e({\a b -c})$ of $\K(\Gbt)$.

We also use the notation of indefinite integrals \footnote{  The useful
  notational element $dx$, along with the conventions of indefinite integration,
led us to adopt integral rather than summation notation.}.   We write:
$$\int f(x) dx = g(x)$$
to mean:  for any $a,b$,  $\int_{a}^{b} f(x) dx = g(b)-g(a)$.

Thus only  $g|^x_y = g(x)-g(y)$ is defined, not $g(x)$.  Nevertheless addition and 
composition on the right with a function make  sense: $(g \circ h) |^x_y = g| ^{h(x)}_ {h(y)}$,
$(g+g')|^x_y = g|^x_y + (g')^x_y$.
  Moreover, if $f,f'$ represent the same element of $Fn(\G,\K(\Gb))$, then for any $h \G \to \G$, 
  $(\int  f(x) dx) \circ h = (\int  f'(x) dx) \circ h$.  
 In particular, $\int_0^x f(t) dt$ induces a well-defined functional 
 $Fn(\G, \K^{df} (\Gb)) \to Fn(\G, \K^{df} (\Gb))$.

\ssec{Dimension filtration}

The ring $\K^{df}(\Gb)$, unlike $\K(\Gb)$, no longer keeps track of ambient dimension;
 but we still have a filtration based on  intrinsic dimension:
 
$$F_n(\K_{\Qq}^{df}(\Gb)) = \{ \a [X] / [0]_1 ^{-m} - \b [Y] /   [0]_1 ^{-m}  :  \a,\b \in \Qq, X,Y \in \Gb[m],
\dim(X),\dim(Y) \leq n \}$$

Let $Gr_n \K_{\Qq}^{df} (\Gb) = F_n(\K_{\Qq}^{df}(\Gb)) / F_{n-1} (\K_{\Qq}^{df}(\Gb))$.

The graded version is not needed at the level of results; but it will simplify the proofs
inasmuch as without it the integration by parts formulas become more complicated.

\<{lem}  
Let $a<b$ be definable points.  There exists a unique linear map
$$gr \int_{a}^b dt:  Fn(\G, Gr_{n-1} \K_{\Qq}^{df} (\Gb)) \to Gr_{n} \K_{\Qq}^{df} (\Gb))$$
such that for any  bounded, definable $X \subseteq \G \times \G^n$, if $\dim X_t \leq n-1$,
and
$f(t)$ is the class of $X_t$ in $Gr_{n-1} \K_{\Qq}^{df} (\Gbt)$,  then 
$$gr \int_{t=a}^b f(t) dt = [X]$$
where $[X]$ is the class of $X$ in $Gr_n \K_{\Qq}^{df}(\Gbt)$. \>{lem}

\ssec{Integration by parts}

 Let $C$ be one of the categories:  $\G_A,\Gb,\volG _A \bdd$.
  Let $\K$ be the Grothendieck ring of $C$.  
 
 The category and the ring $\K$ are then $\Nn$-graded, with a canonical 
homogeneous element $[0]_1$ of grade $1$,  and we can form
the dimension free ring $\K^{df}$.  We also have canonical maps
$\K[n] \to \K[n+1]$, multiplication by $[0]_1$.  Integrals over $\G$ of objects in $\G[n]$ do not
in general exist in $\G[n]$, but if the objects come from $\G[n-1]$ they do; thus 
integration over $\G$ (or a definable interval in $\G$) gives an operator $\G[n-1] \to \G[n]$.   
The integral notation extends formally to $\K^{df}$.

 For $1 \leq i \leq n$, let  $f_i \in Fn(\G,\K^{df}), F_i(x)  = \int_0^{l_i(x)} f_i(t) dt$,
where $l_i$ is a monotone increasing definable function $\G \to \G$.  Also
let $\bF_i(x) = F_i(x) + f_i(l_i(x))$.  

\<{lem} \lbl{ibp-1}  Let $b \in \Qq \tensor A$.  We have equality of classes in $ \K $:
$$ \Pi_i F_i (b) = \sum_{i=1}^n \int_0^{l_i(b)} f_i(t) \prod_{j < i} F_j( l_i \inv (t)) \prod_{j>i} \bF_j( l_i \inv(t)) dt $$
\>{lem}  

\prf  It suffices to prove the same statement for $f_i \in Fn(\G,\K_{+})$, since
it is linear in each $f_i$ and hence formally extends to $\K$, and thence to $\K^{df}$  by
division.  

For $t=(t_1,\ldots,t_n) \in \G^n$, let $i(t)$ be an index $i \in \{1,\ldots, n\}$ with $l_i 
\inv (t_i)$
having the maximal value.  
In case there are several such indices, let $i(t)$ be the smallest possible one.
$\Pi_i F_i$ is the class of
$$ \sum \{f(t_1) \cdot \ldots \cdot f(t_n): t \in X \}$$  
where $X= \{(t_1,\ldots,t_n): 0 \leq t_i < l_i(b) \}$.
Let $X_i = \{(t_1,\ldots,t_n): i(t) = i \}$.  Then $X$ is the disjoint union of the $X_i$, 
and 
$$X_i = \{t: 0 \leq t_i < l_i(b), l_j  \inv (t_j) < l_i \inv (t_i) (j<i), l_j \inv (t_j) \leq l_i \inv (t_i) (j \leq i) \}$$
The formula follows.  \eprf

Now assume in addition that $g \in Fn(\G,\K(\Gb))$, $G(x) = \int_0^{x} g(x)$.

\<{cor} \lbl{ibp-2}  
$$\int_{0}^{b}  g(t) \cdot \Pi_i \bF_i(t) dt  = G(b) \cdot \Pi_i F_i (b)- 
\sum_{j=1}^n  \int_0^{l_j(b)} G(l_j \inv(t)) f_j \Pi_{1 \leq  k < j}  F_k( l_j \inv) \Pi_{j<k\leq n} \bF_k(l_j \inv) dt $$ 

$$\int_{0}^{b}  g(t) \cdot \Pi_i  F_i(t) dt  = G(b) \cdot \Pi_i F_i (b)- 
\sum_{j=1}^n  \int_0^{l_j(b)} \bG(l_j \inv(t)) f_j \Pi_{1 \leq  k < j}  F_k( l_j \inv) \Pi_{j<k\leq n} \bF_k(l_j \inv) dt $$

\>{cor}

\prf  Obtained by subtraction from \lemref{ibp-1} in the case of $n+1$ functions, with $G=F_0$ and $l_0(x)=x$ for the first equation, $G=F_{n+1}, l_{n+1}=x$ for the second. \eprf

We will often look at highest homogenous terms.  The degree will be clear from
the context, so we will write $\eg$ for equality in the graded ring.   In the graded ring 
 there is no distinction between $F_i$, $\bF_i$ and the formula simplifies to:  

\eqn{ibp-3}{ \int_{0}^{b}  g \cdot \Pi_i F_i(t) dt  \eg G(b) \cdot \Pi_i F_i (b) - 
\sum_{j=1}^n  \int_0^{l_j(b)}  G(l_j \inv(t)) f_j \Pi_{1 \leq k \neq j}  F_k( l_j \inv)  dt   }

The variable limits of integration are needed because of the expression below for $\i(\a x+c)$;
it cannot be written as an integral with limits $0,x$ of a function.  

\<{claim} \lbl{dif1} Let $\a = q/p \in \Qq$ be a reduced fraction.  Then 
$$[0, \a x +c) = \int_0^{qx+pc}  e(\frac{x}{p}) dx    $$  \qed \>{claim}

Now in \eqref{ibp-3} we take,   
 for  $i \geq 1$:

 $f_i(x)= e( \frac{x}{p_i} )$ \ ($1 \leq p_i \in \Nn$) 
 
$l_i(x) = q_ix+p_i c_i$ \ ($  c_i \in \Qq \tensor A$, $1 \leq q_i \in \Nn$).  

$\a_i = q_i/p_i$

By  \lemref{dif1} we have $F_i(x) = \int_0^{l_i(x)} f_i(x) = \i( \a_i x + c_i) $.  Hence \eqref{ibp-3} gives:

$${  \int_{0}^{b} g(t) \Pi_{j=1}^n \i(\a_j t+ c_j) dt  \eg G(b) \cdot \Pi_{j=1}^n \i( \a_j b + c_j) - 
\sum_{j=1}^n H_j        }$$

where $H_j =  \int_0^{l_j(b)}  G(l_j \inv(t)) e(t/p_j) \Pi_{1 \leq k \neq j}  F_k( l_j \inv (t))  dt   $.
Now the change of variable $s=t/p_j$ gives:

$$ H_j = \int_0^{\a_j b + c_j} G( \a_j \inv (s - c_j) ) \Pi_{1 \le k \neq j} \i( \a_k (  \a_j \inv(s-c_j)) + c_k )     ds $$

>From this we retain:

\<{lem} \lbl{ibp-5}
$$ \int_{0}^{b} g(t) \Pi_{j=1}^n \i (\a_j t+ c_j) dt  \eg G(b) \cdot \Pi_i \i( \a_i b + c_i) - 
        \sum_{j=1}^n  \int_0^{b_j} G( \frac{s-c_j}{\a_j}) \Pi_{1 \le k \neq j}  \i( \frac{\a_k}{\a_j} s - c_{jk} )ds $$
 
 where  $b_j = \a_j b + c_j$,
$c_{jk} = c_k - \a_j \inv \a_k c_j$.  Note that $c_{jk} = c_k$ if $c_j=0$. \end{lem}

\<{cor} \lbl{powers}   $\int_0^b \i(t)^n dt \eg \frac{\i(b)^{n+1}}{n+1} $ \>{cor}
\prf  Let $g(t)=1, \alpha_j=1, c_j=0$.  Then $G(b)=\i(b)$, and \lemref{ibp-5} gives:

$$\int_0^b \i(t)^n \eg \i(b) \i(b)^n dt -  n \int_0^b \i(s) \i(s)^{n-1} ds = \i(b)^{n+1} - n \i(s)^n $$

Changing sides, we obtain $(n+1) \int_0^b \i(t)^n dt \eg \i(b)^{n+1}$, whence the corollary. \eprf

We will need a more precise version later.   In any $\Qq$-algebra, one can define  
$c_n(x) := {{x}\choose{n}} = \frac{x(x-1)\cdot \ldots \cdot (x-n)}{n!}$.  Note:
\eqn{comb}{c_{n-1}(t)(t-(n-1)) = n c_n(t)}

  Let $C_n(x) = c_n(\i(x))$.  Thus $C_0(x)=1$, $C_1(x)=\i(x)$.
  
\<{lem}   \lbl{powers+}   For $b \in \Qq \tensor A$, $\int_0^b C_n(t) dt = C_{n+1}(b)$. \>{lem}

\prf For $n=0$ this is clear; we proceed by induction.  By \lemref{ibp-2} with $g(x)=1$,
$\bG(x) = x+1$, 
$l_0(x)=l_1(x)=x$, $f_1 = C_{n-1}$, $F_1 = C_n$, we have:
$$ \int_0^b C_n(t) dt = \int_0^b 1 \cdot C_n(t) dt =  \i(b) C_n(b) - \int_0^b (1+t) C_{n-1}(t) dt$$
Now 
$(1+t) C_{n-1}(t) = (t-(n-1))C_{n-1}(t) + n C_{n-1}(t) = nC_n(t) + n C_{n-1}(t)$.  Thus
using the induction hypothesis and \eqref{comb} for $n+1$,
$$(n+1)  \int_0^b C_n(t) dt = \i(b) C_n(b) - n \int_0^b C_{n-1}(t) dt = \i(b) C_n(b) - n C_n(b) =
(n+1)C_{n+1}(b)$$ \eprf

\ssec{Zero-dimensional functions}
Consider   elements of $Fn(\G,\K^{df}(\Gb))$ of the form $e(\a x + \b a) $,
with $\a,\b \in \Qq$, $a \in A$.
By definition, two such terms $e_1,e_2$ are equal iff for all $M \models DOAG_A$ and $c \in M$, the idempotents $e_1(c),e_2(c)$ are equal elements
of $\K^{df}(\G^{\bdd}_{A(c)})$.  According to   \cite{HK} Proposition 9.2, this in turn holds
iff for all subgroups $T$ of $\Qq \tensor A(c)$ containing $A(c)$, $e_1(c) \in T$
iff $e_2(c) \in T$; In other words, iff $A(c,e_1(c)) = A(c,e_2(c))$.  
More generally,

\eqn{cr1} { \Pi_{i=1}^l e(\a_i x + \b_i a_i ) = \Pi_{i=1}^{l'} e(\a'_i x + \b'_i a'_i ) \in Fn(\G,\K^{df}(\Gb)) }
  iff
for any $c \in M \models DOAG_A$, 
$$A(c,\a_1c+\b_1 a_1,\ldots, \a_l c + \b_l a_l) = 
A(c,\a'_1c+\b'_1 a'_1,\ldots, \a'_l c + \b'_l a'_l)$$

As an application, note the equalities, for $m,m'$ relatively prime integers, $k \in \Zz, b \in A$:

\eqn{cr1.1}{    e( \frac{kx+b}{m} ) e(\frac{kx+b}{m'}) = e(\frac{kx+b}{mm'}) }

\eqn{cr1.2}{  e(\frac{kx+b}{m}) = e(\frac{m'(kx+b)}{m}) }

The term ``piecewise'' will refer to partitions of $\G$ into definable points and open intervals,
including all of $\G$ or half-infinite intervals.  
By a {\em constant term} we mean a piecewise constant function, whose values on each
piece are 
 of the form $e(\frac{b}{m})$ with $m \in \Nn$, $b \in A$. 
 By a {\em standard divisibility term} we mean a term $e(\frac{x+b}{m})e(b)$, with   
 $m \in \Nn$, $b \in \Qq \tensor A$.  The integer $m$ is referred to as the denominator.

\<{lem}  \lbl{1}  Any term $e(\a x + \b b) \in Fn(\G,\K^{df}(\Gb))$ is equivalent
to  a product of a  a constant term with a standard divisibility term.  The   denominator  of
the latter is equal to the   denominator of $\a$ as a 
reduced fraction.  
 \>{lem} 

\prf    The term can be written  as $e(mx+nb)/p$, with $b \in A$, $m,n,p \in \Zz, p  \neq 0$.  Write $m=m_1m_2,  p = m_1m_3$, with   $m_2,m_3$   relatively prime.  
As in \eqref{cr1.1}, we have:

\eqn{1.0}{  e(x)e(\frac{mx + nb}{m_1m_3}) = e(x) e(\frac{nb}{m_1}) e(\frac{m_2x + nb/m_1}{m_3}) }

Now since $m_2,m_3$ are relatively prime, there exists $m' \geq 1$ with $m_2m' =1 \mod m_3$.
In particular, $m',m_3$ are relatively prime.
As in \eqref{cr1.2}, 
\eqn{1.0.1}{ e(\frac{nb}{m_1})  e(\frac{m_2x + nb/m_1}{m_3})   = e(\frac{nb}{m_1}) e(\frac{m'(m_2x+nb/m_1)}{m_3}) = 
 e(m'nb/m'm_1) e(\frac{x+m'nb/m_1}{m_3}) }
This is the product of the constant term $e(m'nb/m'm_1)$ with the standard term
$e(m'nb/m_1)  e(\frac{x+m'nb/m_1}{m_3}) $.
Moreover, $\a = m / (m_1m_3) = m_2/m_3$ has denominator $m_3$.
\eprf

\<{lem} \lbl{zero} Any finite product of terms $e(\a x + \b b) \in Fn(\G,\K^{df}(\Gb))$ 
equals a product of one  standard divisibility term and a number of constant terms.  \>{lem}

\prf   
  
Using \lemref{1} and  \eqref{cr1.1} (with $k=1$), 
 it suffices to consider products of terms $e( \frac{x+b}{m} ) e(b)$ with $m$ a prime power,
 $b \in \Qq \tensor A$.  

If $m | m'$, we have, using Criterion \eqref{cr1}:

\eqn{1.2}{   e(b)e(b') e( \frac{x+b}{m}) e( \frac{x+b'}{m'}) = e(b)e(b') e(\frac{x+b'}{m'} ) e(\frac{b-b'}{m})  }

Thus for each prime $p$, it suffices to consider one term $e(\frac{x+b}{p^l} )e(b)$, i.e. the
highest occuring power can be used to reduce the others to constant terms.  So we need
only consider products of terms $ e( \frac{x+b_i}{m_i}) e(b_i)$ with the $m_i$ relatively prime.  

Now if $m_1,\ldots,m_k$ are relatively prime, find integers $l_j$ with $l_j = \delta_{ij} (\mod m_i)$ (Where $\delta_{ij}$ is the Kronecker delta.)  Given $b_1,\ldots,b_k \in A$, let $b^* = \sum l_i b_i$;
then  
\eqn{1.3}{ \Pi_{i=1}^k e(b_i) e(\frac{x+b_i}{m_i}) = \Pi_{i=1}^k e(b_i) e(b^*) e(\frac{x+b^*}{\Pi_{i=1}^k m_i})} 

This finishes the proof.

\eprf 

\<{cor}  Any element of $Fn(\G,F_0 \K^{df}(\Gb))$  
 is equivalent
to  a $\Qq$-linear combination of products of the form of \lemref{zero} \>{cor}

\prf   $F_0 \K_{\Qq}^{df}(\Gb)$ is generated by the classes of definable points $p=(p_1,\ldots,p_n)$.
Each $p_i$ has the form $c_i / m_i$ with $c_i \in A$, and the class $[\{p\}] = e(p_1) \cdot \ldots \cdot e(p_n)$.    
Thus any $f \in Fn(\G,F_0 \K^{df}(\Gb))$  is piecewise of the form of \lemref{zero};
i.e. there exists a partition $I_1 \du \ldots \du I_k$ of $\G$ such that $f|I_j = e_j$,
with $e_j$ a $\Qq$-linear combination of a finite product of terms $e(\a x + \b b)$.  Now the 
characteristic functions of the $I_k$ are also constant terms, and using them it is clear that
$f$ itself is of the stated form.  \eprf

Zero-dimensional terms inside integrals can now be eliminated as follows. 

\<{lem} $e(b) \int    e(\frac{x+b}{m}) h(x) dx = e(b)  (\int h(mx-b) dx)\circ (\frac{x+b}{m})$ \>{lem}

\prf It suffices
to consider standard divisiblity terms $e(\frac{x+b}{m})$, with $b \in A, m \in \Nn$.  The
substitution $y= (x+b)/m$ leads to:
 
\eqn{2}{  e(b) \int_{x=u}^v   e(\frac{x+b}{m}) h(x) dx = e(b)  \int_{y=\frac{u+b}{m}}^{\frac{v+b}{m}} h(my-b) dy    }  \eprf

Note that the analogous formula with rational $m$ would {\em not } be valid; in effect we used
the fact that 
$e(x)e(b)  e(\frac{x+b}{m}) = e(b) e((\frac{x+b}{m}))$.

We note in passing a more direct approach to the   computation of the length of a  segment on lines through the origin; but this method, that ignores the arithmetic of the inhomogeneous part, does
not work for other segments.

\<{lem}  Let $p,q$ be relatively prime integers.  Then there exists $M \in GL_2(\Zz)$
with  
$M \cdot \col{p/q}{1} = \col{1/q}{1} $.%
\>{lem}

\prf  $GL_2(\Zz)$ acts transitively on primitive integer vectors, since they may be completed
to a lattice basis.  Hence some $M \in GL_2(\Zz)$ takes $(p,q)^t$ to $(1,q)^t$.  Thus $M$ takes a planar line
of slope $p/q$ to one of slope $1/q$.   For lines through the origin, the length is now just the
length of a projection. \eprf

\ssec{One-dimensional functions}

\<{lem}   $Fn(\G, F_1 \K_{\Qq}^{df}(\Gb)$ is generated as a $Fn(\G,F_0 \K_{\Qq}^{df}(\Gb))$-module by the terms $\i(\ax + b)$, $\a \in \Qq$, $c \in \Qq \tensor A$. \>{lem}

\prf   A bounded, definable, one-dimensional subset of $\G^n$ is a finite union of points and bounded segments
on lines in $\G^n$, i.e. additive translates of 1-dimensional definable subspaces
$(\a_1,\ldots,\a_n) \G$, with $\a=(\a_1,\ldots,\a_n) \in \Qq^n$.

  We can take $\a$ 
to be a primitive element of $\Zz^n$.  All such elements are $GL_n(\Zz)$-conjugate,
so in fact we can take $\a = (1,0,\ldots,0)$.  In this case the translate has the form
$\G \times \{p\}$, with $p=(p_2,\ldots,p_n)$ a definable point of $\G^{n-1}$.  So the segment has the form
$(a,b) \times \{p\}$, with $a, b \in \Qq \tensor A$.  Hence the class of the segment is
$[(a,b) \times \{p\}] =( \i(b)-\i(a) - e(a) )e(p_2) \cdot \ldots \cdot e(p_n)$.   So $F_1 \K_{\Qq}^{df}(\Gb)$
is generated as an $F_0 \K_{\Qq}^{df}(\Gb)$-module by the elements $\i(b)$, $b \in \Qq \tensor A$.
The lemma  follows. \eprf

For later use, if $\a = p/d$ with $p,d \in \Nn$, and $b \in A$, we will say that $\i(\ax+b)$ admits internal
denominator $d$.  A product of terms, each admitting internal denominator $d$, will 
also be said to admit this denominator.  Note that in general $\i((1/d) a) \neq (1/d) \i(a)$
(even modulo $F_0$.)
 

 %

\ssec{Integration of higher dimensional functions}

Recall the dimension filtration  $(F_n)_A = F_n(\K_{\Qq}^{df}(\Gb))$.  Let $F_n'(\K_{\Qq}^{df}(\Gb))$ be the
$\Qq$-subspace of $F_n(\K_{\Qq}^{df}(\Gb))$ generated by products
of elements of $F_0(\K_{\Qq}^{df}(\Gb))$ with $\leq n$ elements of  $F_1(\K_{\Qq}^{df}(\Gb))$.  
We seek to show (cf. \propref{meas}) that $F_n=F_n'$, i.e. $F_n(\K_{\Qq}^{df}(\Gb))=F_n'(\K_{\Qq}^{df}(\Gb))$ for all $A$
and $n$.  

Let $\FF_n = Fn(\G,F_n' (\K_{\Qq}^{df}(\Gb)))$.
 We will also use an arithmetic refinement:  let $\FF_{n,d}$ be the $\FF_0$-submodule of $\FF_n$ 
generated by $\FF_{n-1}$   along with $n$-fold products of basic one-dimensional terms
with internal denominator dividing $d$, i.e. terms  $\i(\frac{p}{d} x + b)$, $p \in \Nn$, 
$b \in \Qq \tensor A$.

 \<{lem}
  \lbl{I2p} Let  $d,d',p_i \in \Nn$, $c_i,c \in \Qq \tensor A$, $\a_i = p_i / d$, $\g = d/d'$,
  $$f(t) = \Pi_{i=1}^n  \i(\a_i t + c_i) $$
 Then $\int_0^{\g x+c} f(t)dt  \in \FF_{n+1,d'}$  \>{lem}

\prf  We use induction on $d$.  
 Since $\i((\a+1)t+c_i) = \i(t) + \i(\a t + c_i)$ as functions of $t$ in $Fn(\G,\K_{\Qq}^{df}(\Gb))$, 
and using additivity of the integral, we may assume $p_i \leq d$.  Similarly,
$\i(t+c_i) = [ [t,t+c_i) + [0,t) ]/[0]_1 = [ [0,c_i) + [0,t) ] / [0]_1 = \i(c_i) + \i(t)$; so we may
assume that {\em if } $\a_i=1$ then $c_i=0$.

In case $d=1$, we have $p_i= \a_i=1$,  so $c_i=0$ and  $\i(\a_it + c_i)  = \i(t)$. 
By \lemref{powers}, $\int_0^{\frac{x}{d'}+c} \i(t)^n \eg \frac{1}{n+1} \i({\frac{x}{d'}+c})^{n+1}$.
Clearly this expression lies in $ \FF_{n+1,d'}$.   
 
In general,   let $J_1= \{j \leq n: \a_j = 1 \}$, $J_2 = J \m J_1$.  
For $j \in J_1$ we have $c_j=0$.

Using \lemref{ibp-5} with $g=1$, we have:
$$\int_{0}^{\g x+c}   \Pi_{j=1}^n \i (\a_j t+ c_j) dt  \eg \i({\g x + c}) \cdot \Pi_i \i( \a_i (\g x + c) + c_i) -   \sum_{j=1}^n h_j(\a_j (\g x + c) + c_j) $$
  where 
   $$h_j(y)= \int_0^{y} \i( \frac{s-c_j}{\a_j}) \Pi_{1 \le k \neq j}  \i( \frac{\a_k}{\a_j} s - c_{jk} )ds  
   =  \int_0^{y} \i( \frac{d}{p_j}(s-c_j)) \Pi_{1 \le k \neq j}  \i( \frac{p_k}{p_j} s - c_{jk} )ds
   $$

Now if $\a_j=1$ and $c_j = 0$, then    $c_{jk}=c_k$.   Thus (using also $p_j=d$)  each 
of the terms
$h_j(\a_j (\g x + c) + c_j)$ 
is identical with $\int_{0}^{(\g x + c)}   \Pi_{j=1}^n \i (\a_j t+ c_j) dt $.  Moving these terms to the left
we have, with $\nu = |J_1|+1$, $c_j' = c_j / \a_j$:

$$ \nu \int_{0}^{\g x + c}   \Pi_{j=1}^n \i (\a_j t+ c_j) dt  \eg \i(\g x + c) \cdot \Pi_i \i( \a_i (\g x + c) + c_i) -       \sum_{j \in J_2}  h_j(\a_j (\g x + c) + c_j)  $$
 
For $j \in J_2$ we have $p_j <d$, so the induction hypothesis applies.  
Since $ \a_j \g= \frac{p_j}{ d'} $, we have
$h_j(\a_j \g x + \a_j c) \in \FF_{n+1,d'}$.   The remaining terms $\i(\g x + c)$,
$\i(\a_i \g x + \a_i c + c_i)$  clearly have internal denominator $d'$.   This 
concludes the proof of the lemma.  
\eprf

 \<{lem}
  \lbl{I2} Assume $F_n = F_n'$ for all ordered Abelian groups $A$.   Let  $f \in \FF_{n}$.  Then 
  $\int_0^{x}  f(t) dt  \in \FF_{n+1}$.   
        \>{lem}

\prf  It follows from the hypothesis, applied to the structure generated by an element $b$, that 
$f(b) \in F_n'$; it follows by compactness that $f$ itself is a product of $0$- and $1$-dimensional generators.  
By  \lemref{zero}, any product of 
$0$-dimensional generators equals a product of
one standard divisibility term $e(\frac{x+b}{m})e(b')$ and constant terms.  The constants commute
with integration and may be ignored.  So we may assume $f=e(\frac{x+b'}{m})e(b')g_1\cdot \cdot g_n$ with $g_i$ a basic one-dimensional term.  Now with the change of variable $s = \frac{t+b'}{m}$ we have $e(\frac{t+b'}{m})dt = e(s) ds = ds$, i.e.
$${ \int _0^{b}  f(t) dt = 
                  \int_0^{\frac{b+b'}{m}}   g_1(ms-b') \cdot \ldots g_n(ms-b')  ds   }$$

Since $g_i(ms-b)$ is again a basic one-dimensional term,  we may assume: 
$$f(t) = \Pi_{i=1}^n  \i(\a_i t + c_i) $$
in order to show:  $\int_0^x f(t)  \in \FF_{n+1}$.  This follows from \lemref{I2p}.  \eprf

\<{prop}\lbl{I3}  $\K_{\Qq}^{df}(\Gb)$ is generated as a $\Qq$-algebra by the elements
  $e(a), \i(a)$, $a \in \Qq \tensor A$.  \>{prop}
  
\prf  We have seen that $F_0',F_1'$ are contained in the algebra generated by these terms. 
Hence
  it suffices to show that $F_n=F'_n$ for each $n$.  For $n=0,1$ this is true by definition;
  we proceed by induction.  Assume $F_n=F'_n$, and let   $X \subseteq \G^{n'}$ be
  definable and bounded, of dimension $\leq n+1$.  After a finite definable partition we may assume the first
projection has fibers of dimension $\leq n$.  By induction, for any $t$, $[X_t] \in F_n(\K_{\Qq}^{df}(\Gbt))$.
It follows that there exists a  definable partition
 $\G = \union_j I_j$ and $f_j \in Fn(\G,F_n)$ such that for $t \in I_j$, $[X_t]=f_j(t)$.  We may
 take $I_j$ to be an interval $(a_j,b_j)$ $(j \in J_0)$ or a singleton $\{c_j\}$ ($j \in J_1$), or $f_j=0$.  
 Then $X$ is the disjoint union of the pullbacks of the $I_j$; so we may assume
$$ [X] = \sum_{j \in J_0}  \int_{a_j}^{b_j} f_j + \sum_{j \in J_1} f_j(c_j) $$
By \lemref{I2},  $\int_{a_j}^{b_j} f_j = \int_0^{b_j} f_j -  \int_0^{a_j} f_j - f_j(a_j) \in F_{n+1}'(\K_{\Qq}^{df}(\Gb))$.   Thus $[X] \in   F_{n+1}'(\K_{\Qq}^{df}(\Gb))$.  \eprf

In fact we have obtained a somewhat stronger statement.  The semiring
$\SGo(\G)$ was defined below \defref{Gcat}. Let $\SGo(\G)'$ be the subsemiring generated by
the elements
  $e(a), \i(a)$.  Let  
  $\K^{df}(\Gb)'$ be the corresponding ring, and 
$$  \K_{\Qq}^{df} (\Gb)' = \Qq \tensor \K^{df}(\Gb)'$$
Let $\SGo(\G)''$ be the semiring obtained from $\SGo(\G)$ by adding additive inverses
to the elements of $\SGo(\G)'$, and $  \K_{\Qq}^{df} (\Gb)''$ the result of formally dividing by
integers $n >0$.  

  We have natural homomorphisms
$$\K_{\Qq}^{df} (\Gb)' \to   \K_{\Qq}^{df} (\Gb)'' \to  \K_{\Qq}^{df} (\Gb)$$

\<{lem}  The natural homomorphism 
 $\K_{\Qq}^{df} (\Gb)' \to   \K_{\Qq}^{df} (\Gb)$ is an isomorphism   \>{lem}
 
 Explicitly, for any $a \in \SGo(\G)$ there exists $m \in \Nn$ and $b,c \in \SGo(\G)'$
 such that $ma + b =c $ in $\SGo(\G)$.
 \prf  
All our integral equalities are valid in $  \K_{\Qq}^{df} (\Gb)'$.  Hence the proof of \propref{I3} shows
that $\K_{\Qq}^{df} (\Gb)' \to   \K_{\Qq}^{df} (\Gb)''$ is surjective.  Since the same elements are inverted in these semirings, the homomorphism is also injective, hence  bijective, and
$ \K_{\Qq}^{df} (\Gb)''$ is in fact a ring.  Since $ \K_{\Qq}^{df} (\Gb)$ is obtained from 
$\K_{\Qq}^{df} (\Gb)''$ by additively inverting elements, the homomorphism 
$  \K_{\Qq}^{df} (\Gb)'' \to  \K_{\Qq}^{df} (\Gb)$ is also an isomorphism.  \eprf

 \ssec{Subrings and quotients of $\K_{\Qq}^{df}(\Gb)$}

Let $A$ be an ordered Abelian group, 
and let $T_A$ denote the symmetric algebra  
$\Qq \oplus (\Qq \tensor A) \oplus Sym^2 (Q \tensor A) \oplus \ldots$.   If $A=\Zz^n$,
this is a polynomial ring in $n$ variables.

We have a homomorphism $\phi_A: T_A \to \K^{df}(\Gb))$, $a \mapsto \i(a)$.  The image contains
the classes of points of $A$ (all equivalent to $1$) and segments with endpoints in $A$.

\<{lem}\lbl{subring} The natural homomorphism $\phi_A: T_A \to  \K_{\Qq}^{df}(\Gb)$ is injective.

If $A$ is  divisible , $\phi_A$ is an isomorphism. \>{lem}

\prf    We may assume $A$ is finitely generated.   First consider the case $A \subseteq \Qq$. 
So $A \cong \Zz$, and we may take $A = \Zz$.  The symmetric algebra $T_A$ can be identified
with the polynomial ring $\Qq[T]$.   Given a nonzero polynomial $f \in \Qq[T]$, we must
show that $f(\i(1)) \neq 0 \in \K_{\Qq}^{df}(\Gb)$.  Now for any $m$, we have a homomorphism
$$count_m:  \SG[\Gb]   \to \Qq: \ \  [X] \mapsto \#(X \meet ((1/m)\Zz)^n $$
counting points of a bounded definable set $X \subset \G^n$ with coordinates in $(1/m)\Zz$.
This is clearly $GL_n(\Zz)$-invariant, and induces a ring homomorphism
$count_m:   \K_{\Qq}^{df}(\Gb) \to \Qq$.  Composing with $f \mapsto f(\i(1))$ we have a homomorphism $c_m: \Qq[T] \to \Qq$.   Now $c_m(T) = \#(\i([0,1) \meet (1/m)\Zz) = m$.
So $c_m(f) = f(m)$.  Since $\Zz$ is Zariski dense in the affine line, $c_m(f) \neq 0$ for some $m$.
It follows that $f(\i(1)) \neq 0$.

For the general case we will use a statement of 
Van den Dries,   Ealy, and    Marikova. 
 The proof is included in \cite{HK} Proposition 9.10, with $\Rr$ in place of $\Qq$, but this does not matter.

\Claim{}  Let $Q \in \Qq[u_1,\ldots,u_n]$, $B \subset \G^n$ a $DOAG$-definable set, 
and  $Q$ vanishes on $B(\Qq)$, then $Q$ vanishes on $B$.


An element of $T_A$ can be written as $G(a)$, with $G \in \Qq[X_1,\ldots,X_n]$
and $a=(a_1,\ldots,a_n) \in A$.   Suppose $\phi_A(G(a)) = 0$.  This is due
to a finite number of $GL_k(\Zz)$-isomorphisms and $A$-translations between finite unions of products
of the intervals $[0,a_i)$ and points, and possibly some auxiliary intervals and points with endpoints $a'_1,\ldots,a'_{n'}$, that cancel out.  Hence there exists   $DOAG$-definable set 
$B \subseteq \G^{n+n'}$
such that $(a,a') \in B$, and for any ordered Abelian group $A'$, 
 $\phi_{A'}(G(c))=0$ whenever $(c,c') \in B(A')$.
Now suppose in addition that $G(a) \neq 0$.  Then by the Claim, there exist $(c,c') \in \G(\Qq)^{n+n'}$ with $G(c) \neq 0$.  But $\phi_\Qq(G(c)) = 0$.  This contradicts the case $A=\Qq$
proved above.  

 If $A$ is  divisible , the homomorphism $\phi_A$ is surjective.  This follows from
\propref{I3}:  all $e(a)=1$, while $\i(a) = \phi_A(a)$.
\eprf

Denote $\bT_A = \phi_A(T_A)$.   This is always a split subalgebra of $\K_{\Qq}^{df}(\Gb)$, equal
to it if $A$ is  divisible .  
To clarify the full structure, we ask:

\<{question}   
For $n=2$ we have $2 \i(a/2) = [ [0,a/2) \union (a/2,a] ] = \i(a) + 1 - e(a/2)$; so 
$\i(a/2) = (1/2) (\i(a)+1-e(a/2))$.    Is this the first term of a sequence of polynomial relations?  \>{question}

The proof of \lemref{subring} may give the impression that  specializations
 of finitely generated subgroups of $\G$ into $\Zz$, followed by the maps $count_m$, 
 resolve points on $\K_{\Qq}^{df}(\Gb)$ and thus give decisive information.  This is not the 
 case, as the example below shows.   
 
\<{example}   $$\int \int  (e(\frac{s}{2})-1)(e(\frac{t}{2})-1)(e(\frac{s-t}{2})-1) ds dt $$ evaluates to $0$ under any
$count_m$, for any choice of $s,t \in \Zz$, but is not identically $0$.  \>{example}

Let $L_A$ be the field of fractions of $T_{\bA}$, where $\bA = \Qq \tensor A$.

\<{cor} \lbl{field}   There exists a natural homomorphism $\psi_A: \K_{\Qq}^{df}(\Gb) \to L_A$, injective on the image of $T_A$.  The kernel is generated by the relations
$e(a)=1, n \i(\frac{a}{n})  = \i(a)$.   \>{cor}

\prf   If $A$ is  divisible , the
homomorphism $\phi_A$ of \propref{subring} has an inverse $\psi:  \K_{\Qq}^{df}(\Gb) \iso T_A$.  It suffices to view $\psi$ as a homomorphism
into the field of fractions $L_A$ of $T_A$.  

In general, let $\bA = \Qq \tensor A$.   We have a natural surjective homomorphism
 $\nu: \K_{\Qq}^{df}(\Gb)  \to \K_{\Qq}^{df}(\G_{\bA} \bdd)$, $[X] \mapsto [X]$.  Composing with $\psi_{\bA}$
 we obtain a homomorphism $\psi_A: \K_{\Qq}^{df}(\Gb) \to L_A$ where $L_A = L_{\bA}$. 
 Since $\nu \phi_A  = \phi_{\bar A} |T_A$, $\psi_A \phi_A = \psi_{\bA} \phi_{\bar A} | T_A = Id_{T_A}$.  This proves the injectivity on $T_A$.
 
 The relations $e(a/n)=1, n \i(\frac{a}{n})  = \i(a)$ ($a \in A$) are already in the kernel of $\nu$;
both are seen using the translation $x \mapsto x+a/n$.  These relations suffice
(using \propref{I3}) to reduce any element of $\K_{\Qq}^{df}(\Gb)$ to an element of the image
of $T_A$.  By the injectivity on $T_A$  no further relations intervene.
\eprf   

\>{section}

\<{section}{The measured Grothendieck  ring}

We turn to the dimension-free Grothendieck ring of the category
${\volG_A}[*]$  of \defref{Gcat} (3-5).  When possible we   omit $A$ from the notation.

We begin by representing this Grothendieck ring as a ring of functions under convolution.

Recall the  semigroup of definable functions $\G \to \SG(\G[n])$
of  \secref{definable-functions}.     Define  a convolution
product
$$Fn(\G,\SG(\G[n-1])) \times Fn(\G,\SG(\G[m-1])) \to Fn(\G,\SG(\G[n+m-1]))$$
as follows:  
   if $f$ is represented by a definable $F \subseteq \G \times \G^m$,
in the sense that $f(\g) = [F(\g)]$, 
and $g$ by a definable $G \subseteq \G \times \G^n$, let 
 $$f*g(\gamma) = [\{(\a,b,c): \a \in \G, b \in F(\a), c \in G(\g - \a)\}]$$

To distinguish this semiring from the semiring $Fn(\G, \SG(\G))$ with pointwise 
multiplication, we denote it $\Fnc(\G, \SG(\G))$. 
 
Let $\Fnc(\G,\SG(\G))[*] = \oplus_m  \Fnc(\G,\SG(\G))[m]$,
a graded semiring.

$\Fncb(\G,*)$ are the functions with semi-bounded domain and pointwise bounded
range:

\<{notation}  $\Fncb(\G,*) = \{f \in \Fnc(\G, *) :  (\exists \g_0)(\forall \g < \g_0)(f(\g)=0)\} $ \>{notation}

\<{lem} \lbl{gamma-volume}  

(1)    $\SG {\volG}[n] \iso Fn(\G,\SG(\G[n-1]))$

(2)   For $n \geq 1$,  $\SG {\volG}[n] \bdd \iso  \Fncb(\G, \SG(\G[n-1] \bdd)) $  (Functions into the dimension-free
Grothendieck semigroup of two-sided bounded subsets of $\G^{n-1}$, whose support is bounded below.)

(3)  $\SG {\volG}[*] \bdd \iso  \Nn \oplus \oplus_{n \geq 1} \Fncb(\G, \SG(\G[n-1] \bdd)) $ 

    \>{lem}

\prf  (Compare  Lemma 9.12  of \cite{HK}; we include a proof for completeness.) 

Note first that the linear map $(x_1,\ldots,x_n) \mapsto \sum_{i=1}^n x_i$
 is $GL_n(\Zz)$-conjugate to the map $(x_1,\ldots,x_n) \to x_1$.  Therefore $\volG$
 is isomorphic to the category $\volG '$ defined in the same way, except with maps
 preserving the form $x_1$ in place of $ \sum_{i=1}^n x_i$.  Moreover we
 relabel the variables as $(t,x_1,\ldots,x_{n-1})$.   Given $X \subseteq \G^n$, and $t \in \G$,
 let $X_t = \{(x_1,\ldots,x_{n-1}): (t,x_1,\ldots,x_{n-1}) \in X \}$.  

Given a 
semi-bounded definable $X \subseteq \G^n$ , let $\a(X)$ be the definable function:
$t \mapsto [X_t] $.  If $h: X \to Y$ is a $\volG '$-isomorphism, then clearly
$h$ restricts to bijections $h_t: X_t \to Y_t$ which are in fact $\volG_{A(t)} [n-1]$ -isomorphisms.
Hence $\a(X)$ depends only on $[X]$, and a homomorphism
$$\a: \SG {\volG}[n] \to Fn(\G,\SG(\G[n-1]))$$ 
is induced.

Conversely given $F \subseteq \G \times \G^{n-1}$ representing an element of $Fn(\G,\SG(\G[n-1]))$, let $\b(F) = [F]$, the class in ${\volG}[n]$ of the graph of $F$.  If $F,F'$ represent the same element of $Fn(\G,\SG(\G[n-1]))$, then for any $t$, $F_t,F'_t$ are
$\G_{A(t)} [n-1]$ - isomorphic.  The isomorphism is given by a definable
bijection $g_t: F_t \to F'_t$.  By a standard compactness argument (cf. \cite{HK} Lemma 2.3)
we can take $g_t$ definable uniformly in $t$, and define
$g(t,x_1,\ldots,x_{n-1}) = (t,g_t(x_1,\ldots,x_{n-1}))$; then $g: F \to F'$
is a definable bijection.  Moreover for any $t$, there is a finite set of matrices
$M_1(t),\ldots,M_{k(t)}(t) \in GL_{n-1}(\Zz)$ and elements $c_i(t) \in A(t)$ such that for any $x \in \G^{n-1}$,
for some $i \leq k(t)$, $g_t(x) = M_i(t) x + c_i(t)$.  By compactness, $M_1(t),\ldots,M_k(t)$
can be chosen from a finite set $M_1,\ldots,M_k$ of matrices.  So for any
$t \in \G$ and $x \in \G^{n-1}$, for some $i \leq k$, $g_t(x) - M_i(t) x \in A(t)$.
Now $A(t) $ is the group generated by $t$ over $A$, so any element of $A(t)$ has
the form $a+mt$ for some $m \in \Zz$. By compactness, there exist finite subset
$A_0 $ of $A$ and $Z_0$ of $\Zz$ such that 
 for any
$t \in \G$ and $x \in \G^{n-1}$, for some $i \leq k$, some $a \in A_0$  and $m \in Z_0$, 
$g_t(x) =  M_i(t) x + a+ mt$.  Partition $X$ into finitely many pieces, such that $M_i, m,a$
are constant on each piece; then on each piece $g$ is given by 
$(t,x) \mapsto (t,M x + a + mt)$for some $a \in A^{n-1} $ and $m \in \Zz^{n-1}$.
But this is clearly an affine $GL_n(\Zz)$-transformation.  Thus $g$ is a $\volG ' [n]$-isomorphism.
So $[F]=[F']$ in $\volG ' [n]$.   his allows us to define $\b: Fn(\G,\SG(\G[n-1])) \to \SG \volG ' [n]$.
 
It is clear that $\a,\b$ are inverse homomorphisms.  So $\a$ is an isomorphism and shows (1).

Restricting $\a$ to bounded sets yields an isomorphism yields (2).
$$\SG {\volG} [n] \bdd \to  \{f \in Fn(\G, \SG(\G[n-1] \bdd)):  (\exists \g_0)(\forall \g < \g_0)(f(\g)=0)\} $$
The direct sum of these isomorphisms over all $n \geq 1$  gives (3);  in grade $0$ we have $\Nn$ on both sides.  The verification that product
goes to convolution product is straightforward.     \eprf

 Let $q_0 \in  \Fncb(\G, \SG(\G[0] \bdd))$ be the function with support at $\{0\}$
and value $1$.   Note that for $f \in \Fncb(\G, \SG(\G[n] \bdd))$,
 $f*q_0$ is the element of $ \Fncb(\G, \SG(\G[n] \bdd))$ satisfying
 $(f*q_0)(t) = f(t) \times [\{0\}]_1$.

\def\Fn{{\rm Fn}}
We can also define a convolution product on the semigroup
 $ \Fn \bdd (\G, \K_{+}^{df}(\G \bdd))$.
An element of this semigroup is represented by a pair $(f,n)$, where
$f \in \Fn \bdd (\G,\K_{+}(\G[n] \bdd))$, and 
$(f,n)$ is identified with $(f * q_0^m,n+m)$.  The pair $(f,n)$ is  intended to represent the function
$t \mapsto f(t) [0]_1 ^{-n}$.     
We let $(f,n)*(g,m) = (f*g, n+m+1)$.  This makes  $ \Fn(\G, \K_{+}^{df}(\G \bdd))$
into a semiring $ \Fnc(\G, \K_{+}^{df}(\G \bdd))$.

\<{lem}  \lbl{gamma-df}    $\K_{+}^{df} (\vol \G \bdd)$ is canonically
isomorphic to   $  \Fncb(\G, \K_{+}^{df}(\G \bdd))$
\>{lem}

\prf     
By \lemref{gamma-volume} (3), 
$$\K_{+}^{df} (\vol \Gb)[*] \cong (\oplus_n \Fncb(\G, \SG(\G[n-1] \bdd)))[q_0 \inv]_0$$
Let $(f,n)$ represent an element of $  \Fncb(\G, \K_{+}^{df}(\G \bdd))$. Let
$$ [(f,n)] \mapsto f q_0^{-(n+1)}$$
This defines an injective semiring homomorphism
$$  \Fncb(\G, \K_{+}^{df}(\G \bdd)) \to  (\oplus_n \Fncb(\G, \SG(\G[n-1] \bdd)))[q_0 \inv]_0$$
which is clearly also surjective.

\eprf

Given a definable function $h: \G \to \Fnc(\G, \SG(\G))$, and a definable $Y \subseteq \G$,
we define $\int_Y  h \in \Fnc(\G, \SG(\G))$ pointwise, i.e.  $(\int_Y h) (\g) = \int_{t \in Y} ev_\g(h) (dt)$
where $ev_\g(h)(t) = h(t)(\g)$.  
 This carries over to the groups and rings considered below.

Let $\RG= \Qq \tensor \Fncb(\G, \K_{+}^{df}(\G \bdd)) = \Fncb (\G, \K_{\Qq}^{df}(\G \bdd))$ be 
$\Qq$-algebra of functions represented by elements whose support is bounded below.
Then $\RG$ also has a natural convolution structure, and forms a ring.   We begin by developing some identities in $\RG$.   We denote convolution of functions $f,g$ by
$fg$; we will not consider the pointwise product except when one of the functions
is supported on $\{0\}$, in which case the two products are equal.

Let $\RG_0$ be the subring of $\RG$ consisting of elements with support $ \{0\}$ (and
$0$.)
The map $a \mapsto a q(0)$   gives an homomorphism of rings
$\K_{\Qq}^{df}(\Gb) \to  \RG_0$.    In fact, since equality of functions in $\Fn(\G, \K_{+}^{df}(\G \bdd)) $ is defined pointwise, and implies equality of the value at $0$, it is easy to see that
this is an isomorphism.  

\eqn{meas0}{   \K_{\Qq}^{df}(\Gb) \cong  \RG_0 }

Let $q(\g)$ denote the element supported
on $\{\g\}$, with $q(\g)=1$.   

  Then $e(\g) q(\g+\g') =q(\g)q(\g')$.    We have
\eqn{I2.5}{ f=\int_{t \in \G} f(t) q(t) dt }
 
$$\int f(t) q(m t) dt = \int f(t/m) e(t/m) q(t) dt$$

The elements
of $\K_{\Qq}^{df}(\G \bdd)$ can  are identified with constant functions with support $\{0\}$.

For $m \geq 1$, and $b \in \Qq \tensor A$, let
  
$\theta_{m,b} = \int_{t \geq b}   q( m t) dt$   
and $\theta_m = \theta_{m,0}$, 
$\theta=\theta_1$. 
Let $Q_m(b) = \int_0^b q(mt) dt$.  So 
$Q_m(b) = \theta_m - \theta_{m,b}$.

The filtration on $ \K_{\Qq}^{df}(\G \bdd))$ induces a filtration $F_n\RG$ on $\RG$.  
$F_0\RG$ consists of ``purely exponential'' sums; it has as a $\Qq$- basis the
elements $\theta_m$, $q(b)$, $Q_m(b)$.  
  Let $F_n\RG_0 = F_n\RG \meet \RG_0$.  
 
Let $F_n'\RG$ be the
$\Qq$-space  generated by    
  products  $q(b') a_1 \cdot \ldots \cdot a_{n}$, 
 where $a_i \in F_1\RG_0$ or $a_i=\theta_{m,b}$ for some $m$ and some 
$b \in \Qq \tensor A$.
 
As above we will write some of 
the identities in graded form.

Note that $e(b) \theta_{m,b} = e(b) \int_{b}^\infty q(mt) dt = e(b) \int_{mb}^\infty e(\frac{s}{m}) q(s) ds $.
Since $e(b) e(\frac{s+mb}{m}) = e(b) e(\frac{s}{m})$, we have: 
$e(b) \theta_{m,b} = e(b) \int_0^\infty e(\frac{s+mb}{m}) q(s+mb) ds  
= e(b) q(mb) \int_0^\infty e(\frac{s}{m}) q(s) ds = e(b)q(mb) \theta_m$.  Hence
 
\eqn{meas1}{   e(b)Q_m(b) =  e(b) \int_0^b q(mt) dt = e(b)  (1-q(mb))  \theta_m } 

Note that while $\int_0^\infty q(t) f(t)$ is defined, $\int_0^\infty f(t)$ is not.  Thus integration
by parts does not directly apply.  To compute unbounded integrals (when $A \neq (0)$)
we will use:

\<{lem} \lbl{meas2}    Let $f(x) = \Pi_{i=1}^n \i(\a_i x + c_i)$.  Let $m \in \Nn$ be such that
$m \a_i \in \Nn$, and let $a \in mA, a \neq 0$.  Then $f(t-a) = f(t) -f_1(t)$ for some
$f_1 \in F_{n-1} \K_{\Qq}^{df}(\Gbt)$;
and we have:
   
$$(1-q(a)) \int_0^{\infty} f(t) q(t) dt = 
\int_0^\infty f_1(t) q(t) dt - \int_0^a f_1(t)q(t) dt + \int_0^a f(t)q(t) dt $$

\>{lem}

\prf  Let $\a_i = p_i / m$, $a = mb$.  We have $\i(\a_i(x-a) + c_i) = \i(\a_ix + c_i -p_ib) = \i(\a_ix-c_i) - p_i \i(b)$.  From this the existence of $f_1$ is clear.  We compute:
$$ q(a) \int_0^{\infty} f(t) q(t) dt  = \int_0^\infty f(t) q(t+a) dt = \int_{a}^\infty f(s-a) q(s) ds = $$
$$= \int_{a}^\infty f(s) q(s) ds - \int_{a}^\infty f_1(s) q(s) ds = \int_0^{\infty} f(t) q(t) dt  - \int_0^{a} f(t) q(t) dt + \int_0^a f_1(t) q(t) dt + \int_0^a f_1(t) q(t) dt$$ 
and the lemma follows.
\eprf

Assuming $A \neq 0$, fix an element $0< a_0 \in A$.  
Define
$\RG_{l} = \RG[(1-q(ma_0))^{-1}: m \in \Nn ]$.  Since the elements inverted are from $F_0\RG$,
the filtration carries through to $\RG_{l}$.  
 Let $\RGb$ be the subring of $\RG$ consisting of elements 
with two-sided bounded support:
$$\RGb = \{\int_{-b}^b f(t) q(t) d: f \in \RG, b \in \Qq \tensor A \} $$
and $\RG_{b,l}$ be the localization of $\RGb$ obtained by inverting the elements
 $(1-q(ma_0))$, $m \in \Nn$.

\<{cor} \lbl{bdd} Assume $A \neq (0)$.  
 Then the inclusion  $\RGb \to \RG$ induces
an isomorphism 
  $\RG_{b,l} \to \RG_{l}$.

  \>{cor}

\prf The surjectivity is clear from \lemref{meas2}, and induction.   Moreover, inspection of the proof shows that additive inverses are used
only for elements in the image of $\RG_{l}$.  Thus if $A$ is the subsemiring of $\RG_{l}$, 
generated by $\Fncb(\G, \K_{+}^{df}(\G \bdd))$ and by the image of $\RG_{b,l}$, then the same induction shows that 
$\RG_{b,l} \to A$ is surjective; and in this case injectivity is evident.     Hence $\RG_{b,l} \cong A = \RG_{b,l}$.  \eprf

Now an analog of \lemref{ibp-5}.   We use integration by parts in $\K(\vol \G_A)$.  Products refer to the Grothendieck
ring of these categories, or equivalently to convolution from the point of view of \lemref{gamma-volume}.

\<{lem} \lbl{ibp-6} Let $\a=\a_0, \a_1,\ldots, \a_n \in \Qq^{>0}, c=c_0,c_1,\ldots,c_n \in \Qq \tensor A$,
 $b_j = \a_j b + c_j$, $c_{jk} = c_k - \a_j \inv \a_k c_j$.   Then 
$$  \int_{0}^{b}  Q_m(\a t +c) \Pi_{j=1}^n \i (\a_j t+ c_j) dt  \eg  
 \i(b) Q_m(\a b+c) \Pi_{i=1}^n \i( \a_i b + c_i)  +$$
 $$ - \int_0^{b_0} q(ms)  \i( \frac{s-c_0}{\a_0})  \Pi_{j=1}^n    \i( \frac{\a_j}{\a_0} s - c_{0k} )ds            +$$
 $$           \sum_{j=1}^n  \int_0^{b_j} \i( \frac{s-c_j}{\a_j}) Q_m(\frac{\a}{\a_j} s - c_{j0}) \Pi_{1 \le k \neq j}  \i( \frac{\a_k}{\a_j} s - c_{jk} )ds  $$
 
  \>{lem}
  

\prf  
Let $F_0(t)=Q(\a t + c)$, $l_0(t)={\a t + c}$,
$f_0(t) = q(mt)$.  
We have by definition 
$F_0(t)  = \int_0^{\a t + c} q(ms) ds$.  We apply \eqref{ibp-3} (for indices $0,\ldots,n$)  with $g=1$, $G=\i$, $f_0,F_1$ as above, and for $i \geq 1$, writing $\a_i = q_i/p_i$,
 $f_i(x)= e( \frac{x}{p_i} )$  
$l_i(x) = q_ix+p_i c_i$, $F_i(x) = \i(\a_ix + c_i)$ as in the proof  of \lemref{ibp-5}.   Thus:
 
$${  \int_{0}^{b} Q_m(\a t +c)  \Pi_{j=1}^n \i_i(\a_j t+ c_j) dt  \eg \i(b)  Q_m(\a b+c) \Pi_i \i( \a_i b + c_i) - H_0 -
\sum_{j=1}^n H_j        }$$

where
$$H_0 = \int_0^{\a b + c}     \i(l_0 \inv (t)) q(mt) \Pi_{k=1}^n F_k(l_0 \inv(t)) dt                        $$
while  for $j \geq 1$,  $H_j =  \int_0^{l_j(b)}  \i(l_j \inv(t)) e(t/p_i)Q_m(\a l_j \inv (t) + c) \Pi_{1 \leq k \neq j}  F_k( l_j \inv (t))  dt   $.
Now $Q_m(\a l_j \inv (p_j s) + c) = Q_m( \a  \a_j \inv (s-c)$, so
the change of variable $s=t/p_j$ gives:

$$ H_j = \int_0^{b_j} \i( \a_j \inv (s - c_j) )Q_m( \a  \a_j \inv (s-c) + c) \Pi_{1 \le k \neq j} \i( \a_k (  \a_j \inv(s-c_j)) + c_k )     ds $$

\eprf

  \<{lem} \lbl{meas-p}    $b,c_i,c \in \Qq \tensor A$, $\a,\a_i \in \Qq^{>0}$, 
 $b \in \Qq \tensor A$. Then 

(1)  $\int_0^{b} q(mt) \Pi_{i=1}^n  \i(\a_i t + c_i)   dt \in F_n' \RG$.

(2)    $\int_0^{b} Q_m(\a t +c) \Pi_{i=1}^{n-1}  \i(\a_i t + c_i)   dt  \in F_{n}' \RG$.  \>{lem}

\prf   
 Let  $d,p,p_i \in \Nn$, $\a_i =    p_i / d$, $\a = p/d$.  
     We use induction on $n$ and on $d$.

If $M_i$ is the largest integer $\leq \a_i$, we have:
   $\i(\a_it + c_i) = M_i \i(t) + \i((\a_i - m)t + c_i)$.   Using this relation, 
   we immediately reduce to the case $p_i \leq d$.    
 
(2.1)  We begin with (2) in the case:  $\a=1$.

We  have
$$Q_m(t+c) = \int_0^{t+c} q(ms) ds =   Q_m(t)+  \int_{t}^{t+c} q(ms) ds $$
Now $e(t) e(mt) = e(t)$, and $e(mt) q(m(t+s)) = q(mt) q(ms)$.  Thus

$$e(t) \int_t^{t+c} q(ms) ds = e(t) \int_0^c q(m(t+s)) ds = e(t) q(mt) \int_0^c q(ms) ds =
e(t) q(mt) Q_m(c)$$
Recall  \eqref{meas1}:
$e(x) Q_m(x) = e(x) (1-q(mx)) \theta_m$.
So
$$ {Q_m(t+c) e(t) = Q_m(t) e(t) +  e(t)q(mt)Q_m(c) = e(t)(1-q(mt))\theta_m + e(t) q(mt)Q_m(c)  }$$
 
Thus the integral (2) equals:

 $$ \theta_m \int_0^{b} (1-q(mt))   \Pi_{k=1}^{n-1}  \i( {\a_k} t - c'_{k} ) dt
 -  Q_m(-c) \int_0^{b} q(mt)  \Pi_{k=1}^{n-1}  \i( {\a_k} s - c'_{k} ) dt$$

Both summands lie in ${F_n' \RG}$, by induction on $n$, and using \propref{I3}.  This finishes
(2) in the case $\a=1$.

(1)     
Let $b_j  = \a_jb + c_j$,
$f(t) = \Pi_{i=1}^n  \i(\a_i t + c_i) $.
    By
\lemref{ibp-5} with $g(t) = q(mt)$, 
$$\int_0^b  f(t) q(mt) e(t) dt     \eg Q_m(b) \cdot \Pi_i \i( \a_i b + c_i) - 
        \sum_{j=1}^n     \int_0^{b_j} Q_m( \frac{s-c_j}{\a_j}) \Pi_{1 \le k \neq j}  \i( \frac{\a_k}{\a_j} s - c_{jk} )ds $$
  The first summand on the right is evidently in $F_n' \RG$.  If $\a_j=1$,  so is the second,
  by the
  case (2.1).  If $\a_j < 1$, then $\a_k / \a_j = p_k/p_j$ have denominators
  $<d$, so induction on $d$ applies and (2) can be quoted.   
     Hence  
 $\int_0^b  f(t) q(mt) e(t) dt \in F_n\RG$.

(2)  in the general case.  We use \lemref{ibp-6} (for $n-1$).  The first summand
on the right is clearly in $F_n' \RG$.   By (1), so is the second.  The remaining $n-1$ summands
are 
$$E_j = \int_0^{b_j} Q_m(\frac{\a}{\a_j}s - c_{j0}) \i(\frac{s-c_j}{\a_j}) \Pi_{1 \leq k \neq j} \i(\frac{\a_k}{\a_j} s - c_{jk}) ds $$
If $\a_j \neq 1$ then again the denominators are $<d$, and by induction 
$E_j \in F_n' \RG$.  If $\a_j =1$ then $E_j$ has the form (2), and so can be moved to the left
as in \lemref{I2p}.

   \eprf

\<{prop} \lbl{meas}  
 
 (1)   $\RGb[\theta_m: m=1,2,\ldots]$ is generated as a $\RG_0=K_\Qq(\Gb)$-algebra by the elements $q(b)$, $\theta_m$  and $Q_m(b)$, $m \in \Nn$, $b \in \Qq \tensor A$. 

 (2)    If $A \neq (0)$,  $\RG_{l}$ is generated  over $\RG_0=K_\Qq(\Gb)$ by the elements
 $q(b), (1-q(ma_0)) \inv$ and
 $\theta_{m,b}$.
  \>{prop}
  
\prf  (1) 
  By \corref{bdd}, an element of $\RGb$ can be written as  a difference of elements
  $\int_0^b  f(t) q(t) dt$, with
$f \in \Fnc(\G, \K_{\Qq}^{df}(\Gb))$.  By the proof of \propref{I3} we can take $f$ to be a product of
zero-dimensional terms and basic one-dimensional terms, restricted to a point or an interval.
Multiplying by an appropriate $q(c)$, we may assume the point is $0$, or the interval
is of the form $[0,c)$.  If $b \leq c$ then the interval may be ignored; if $c < b$
we replace the integral by   $\int_0^b  f(t) q(t) dt$, so that $f$ is defined on $[0,b)$.
    Moreover the 0-dimensional
terms can be collected together to form one basic term $e(\frac{t+b}{m})e(b)$.  
But
$$ \int e(\frac{t+b}{m})e(b) f(t) q(t) dt = \int e(\frac{s}{m})f(s) q(s) ds = \int f(t) q(mt) e(t) dt$$
So it suffices to show that  for 
  $\a_i \in \Qq$, $c_i \in \Qq \tensor A$,
  $$\int_0^b \Pi_{i=1}^n \i(\a_it+c_i) q(mt) e(t) dt$$
 lies in the $\Qq$-algebra   generated
 by the elements $q(b)$  and
 $\theta_{m,b}$, $m \in \Nn$, $b \in \Qq \tensor A$.  This follows from \lemref{meas-p}.  
 
 (2)  Follows from \lemref{meas2}.   \eprf

The next lemma suggests a way to look at unbounded functions; it will not be used further on.
Let $\RGi = \Qq \tensor \Fn(\G, \K_{+}^{df}(\G \bdd)) \cong \Fn(\G, \K_{\Qq}^{df}(\G \bdd))$.
   $\RGi$ is  a $\K_{\Qq}^{df}(\G \bdd)$-module, under pointwise multplication, and more generally
   an $\RG$-module, under convolution.  
Thus we can define 
 $\RGi_{l} := \RGi[ (1-q(ma_0))^{-1}: m \in \Nn ]$.   Note that $\RG$ is no a priori a ring.  However,
 it can be made into one using:
     
\<{lem} \lbl{meas+} Let $0 \neq a \in A$.  The natural inclusion   $\RG_{l} \to \RGi_{l}$ is an $\RG$-module isomorphism.
\>{lem}  

\prf  Using   \propref{meas} together with the automorphism $\g \mapsto -\g$,
 the elements with negative support are generated by the $q(\g)$ together with the elements
  $\theta^-_{m,b}:= \sum_{\g < b}   e(\frac{\g}{m}) q(\g)$.
    But $\theta^-_{m,b}   + \theta_{m,b} =  \sum_\g   e(\frac{\g}{m}) q(\g)  = 0$.
 Since $q(ma) \sum_\g   e(\frac{\g}{m}) q(\g) =  \sum_\g   e(\frac{\g}{m}) q(\g +ma) 
 = \sum_\g e(\frac{\g+ma }{m}) q(\g +ma) = \sum_\g (\frac{\g}{m}) q(\g) $, we have
 $(1-q(ma)) ( \sum_\g   e(\frac{\g}{m}) q(\g) ) = 0$, so in the localized ring we have
$\theta^-_{m,b}   + \theta_{m,b} =0$.  Thus $ \theta^-_{m,b}  = -\theta_{m,b}$ lies
in the image of the localization of $\RG$.  
\eprf

\ssec{The elements $\th{m}$}

We will  see later (\lemref{theta-trans}) that $\theta$ is transcendental over the elements $q(a)$ and $i(a)$.  But
the  various $\theta_m$ are rational over $\theta$:


\<{lem}\lbl{theta} The identity  
$\th{n+1} (\th{} + \th{n} -1) =   \th{n} \th{}$ is valid in $\RG$.  Hence $\th{n}$ is invertible in  in $\RG[\th{} \inv]$, and we have

$$(1- \th{n} \inv)  =    (1- \th{} \inv)^n$$
\>{lem}

\prf  
We have:
\begin{eqnarray} \lbl{t2.1}
\th{n} \th{n+1} & = & \th{n} + \int_{t=0}^\infty q(t) \int_0^t e(\frac{s}{n})e(\frac{t-s}{n+1}) ds dt \\
 \lbl{t2.2}
  \th{n} \th{} & = &  \th{n} +  \int_{t=0}^\infty q(t)  \int_0^t  e(\frac{s}{n}) ds  =  \th{n}+\int_{t=0}^\infty q(t) \i(\frac{t}{n}) dt    \\
   \lbl{t2.3}
  \th{n+1}\th{} & = &   \th{n+1} + \int_{t=0}^\infty q(t)  \i(\frac{t}{n+1}) dt  
\end{eqnarray}
Now by  \eqref{1.2},
\eq{
e(s)e(t) e(\frac{s}{n}) e(\frac{s-t}{n+1}) = e(s) e(t) e(\frac{s+nt}{n(n+1)})
}
With the change of variables $s'=s+nt$ we obtain $e(s)e(t) = e(s')e(t)$, and 
\eq{
\int_0^t e(t) e(\frac{s}{n})e(\frac{t-s}{n+1}) ds  = e(t)  \int_{nt}^{(n+1)t} e(\frac{s}{n(n+1)}) ds
}
With a further change of variable $s'' = \frac{s}{n(n+1)}$, 
\eq{
\int_0^t e(t) e(\frac{s}{n})e(\frac{t-s}{n+1}) ds  = e(t)  \int_{\frac{t}{n+1}}^{\frac{t}{n}} e(s) ds 
= \i(\frac{t}{n}) -   \i(\frac{t}{n+1})
}
so by \eqref{t2.1},
\eq{
\th{n+1}\th{n}   =  \th{n} + \int_{t=0}^\infty q(t) [ \i({\frac{t}{n}}) -   \i(\frac{t}{n+1}) ] dt 
}
  
By    \eqref{t2.2}, \eqref{t2.3}, 
\eq{ 
\th{n+1}\th{n} - \th{n} \th{}  + \th{n+1} \th{} =   \th{n+1}
}
This identity is equivalent to the one in the statement of the lemma.  From this we see that
 $\th{n+1} \inv \in \Qq[\th{},\th{n}, \th{} \inv, \th{n} \inv]$
and 
\eq{
1- \th{n+1} \inv   =    (1 - \th{n} \inv) (1 - \th{} \inv)
}
The lemma follows by induction.   \eprf

\<{remark} \label{localization} \rm  \lemref{theta} will be used to show, after an appropriate localization, that elements of finite support
generate the entire value ring of $\Gamma$.  This will go over to $\RV$ and to $\VF$.  It is a generalization of the {\em rationality of 
Poincar\'e series} and similar rationality results for generating series, in the integration theory of Denef, Denef-Loeser, and Cluckers-Loeser
(polynomials are power series with bounded support.)    \>{remark}

\<{remark} \label{localization2}  \rm
 Another aspect of    \lemref{theta}
  is that a certain localization of the elements of bounded support is forced geometrically.
 The element $\theta_m$
can be written, according to the lemma, as $ \frac{\theta^m}{\theta^m - (\theta-1)^m} $, or again as $\frac{(\G^{\geq 0})^m}{(\G^{\geq 0 })^m - (\Gamma^{>0})^m} $.  

Now (say over an algebraically closed field) the value ring of $\VF$ can be obtained as a tensor product of value rings of $\RES$ and of $\G$, modulo two linear homogeneous relations of degree one:  the equality of the point $0 \in \Gamma$ with the class of the variety  $[G_m]$ over the residue field, 
and the equality of the point $1_k$ of the residue field with the class $[\Gamma^{>0}]$.    Using these relations, we find in $K(\vol VF \bdd)[*][[G_m]_1 \inv ]$ an element whose volume is:
$\frac{[\Aa^1]^m}{[(\Aa^1)^m - [1_k]^m]}$, or $\frac{L^m}{L^m-1}$ in common notation.    Subtracting $1$ we find an inverse of $L^m-1$.  It follows that for any cyclotomic polynomial $c_m(z) = (z-1) \inv (z^m-1)$ with $m >1$,
$h_m([\Aa_1]_1)$ is invertible in  $K_\Qq(\vol VF \bdd)[*]$; though it is not invertible in $K_\Qq(Var_k)$.
%
 \>{remark}

\ssec{Unbounded sets}

We briefly pause to describe the dimension-free Grothendieck ring of $\G$.  
The resulting homomorphisms on $\K(VF)$ were already described in \cite{HK};
the present results confirms their uniqueness.   Compare \cite{marikova}, \cite{kageyama-fujita}.

We denote $e(a)=[\{a\}_1]/ [\{0\}]_1, \i(a) =  [0,a)_1 / [0]_1 $, $\i(\infty)= [0,\infty)_1 / [0]_1 $.

\<{thm} \lbl{unbdd}  $\K_\Qq^{df}(\G_A)$ is generated as a $\Qq$-algebra by the elements $e(a),\i(a)$ 
  \,  $(a \in \Qq \tensor A)$ and $\i(\infty) $.
  
For $a \in A$, we have $\i(a)=0$.   Also $\i(\infty)^2 = -\i(\infty)$.  

If $A$ is  divisible , then $\K_\Qq^{df}(\G_A) \cong \Qq^2$.  \>{thm}

\prf  The proof of  \lemref{ibp-1} remains valid for 
$\K(\G)$ with $b=\infty$, letting $F_i(\infty)=\bF_i(\infty) = \int_0^\infty f_i(t) dt $; and 
the subsequent lemmas through \propref{I3} go through verbatim.  
 This shows that $\K_\Qq^{df}(\G_A)$  is generated
by the elements $e(a), \i(a) $ and $\i(\infty)$.

 The translation $x \mapsto x+a$ shows that   $[0,\infty) = [a,\infty)$.  Hence
 $[0,a)=0$ in $\K(\G[1])$, so $\i(a)=0$.    See \cite{HK} Proposition 9.4 for the relation $\i(\infty)=0$. 
 
 Thus if $A = \Qq \tensor A$, $\K_\Qq^{df}(\G_A)$ is generated by the element
 $\i(\infty)$.  The relation $\i(\infty)^2=\i(\infty)$ shows that the $\Qq$-algebra is a quotient of  $\Qq^2$; the two Euler characteristics   in \cite{HK} show that it is in fact $\Qq^2$.   \eprf
 
\ssec{Subrings and quotients of $\K_{\Qq}^{df} (\vol \G \bdd)$}

Recall \lemref{gamma-df}:  $\K_{\Qq}^{df}(\volG \bdd) := \Qq \tensor  \K_{+}^{df}(\volG \bdd)
 = \Fncb (\G, \K_{\Qq}^{df}(\G \bdd)) =: \RG$.

 Let $L_A$ be the field of \corref{field}.  Let $\bA = \Qq \tensor A$, and let $L_A[q^{\bA}]$ be the formal Puiseux polynomial ring over $L_A$ (i.e. the group ring of $(\bA,+)$ over $L_A$).
Let  $L_A(q^{\bA})$ be the field of fractions.   
Also form the polynomial ring $L_A(q^{\bA})[\theta]$, and rational function field  $L_A(q^{\bA})(\theta)$ (with $\theta$ viewed as an indeterminate.)

\<{lem} \lbl{theta-trans}
$\theta$ is transcendental over $L_A(q^{\bA})$.
 \>{lem}

\prf $\theta^n = \int_0^\infty j(t) q(t) dt$ with $j$ of degree $n$.  Convolving by an
element of $L_A(q^{\bA})(\theta)$ still leaves an expression of the same form, with $j(t) \in F_n \setminus F_{n-1}$.  
The lemma follows 
from the linear independence of polynomials of distinct degrees over the functions
with finite support.  \eprf

\<{prop}\lbl{psi}  Assume $A \neq 0$.  There is a natural homomorphism 
$$\psi_A^*: \RG=  \K_{\Qq}^{df}(\volG \bdd)
 \to  L_A(q^{\bA})(\theta)$$
as well as a homomorphism $\psi_A:  T_A[q^A][\theta] \to \RG$, with
$\psi_A^* \psi_A = Id$.  

  If $A$ is  divisible , $\psi_A^*$ induces an isomorphism
  $$ \RG_{l}[\theta \inv] \to      T_A[q^A][\theta, \theta \inv, (1-q(ma_0))^{-1}, (1-(1-\theta \inv)^m) \inv ]_{m=1,2,\ldots}       $$
    \>{prop}
\prf
Composing the map $\phi_A: T_A \to \K_{\Qq}^{df}(\Gb)$ with the homomorphism 
$\K_{\Qq}^{df}(\Gb)  \to \RG_0$ of \eqref{meas0}, we obtain a map $\psi_A: T_A \to \RG_0$.
We have $\RG =  \Fncb (\G, \K_{\Qq}^{df}(\G \bdd))$.
Extend $\psi_A$  to a homomorphism $\psi_A: T_A[q^A] \to \RG$ with $q^a \mapsto q(a)$.
It is clear by support considerations, and using \lemref{subring}, that $\psi_A$ is injective
on $T_A[q^A]$.  Extend $\psi_A$ further to the polynomial ring $T_A[q^A][\theta]$
mapping $\theta \to \theta$.  By \lemref{theta-trans}, $\psi_A$ remains injective.

Next using \lemref{theta}, extend $\psi_A$ to 
$$\psi_A': T_A[q^A][\theta, \theta \inv, (1-(1-\theta \inv)^n) \inv ]_{n=1,2,\ldots} \to \RG[\theta \inv]$$
  It is still injective, by \lemref{theta-trans}.  
By \lemref{theta}, the image of $\psi'$ contains $\theta_n$ for each $n$.  
By \eqref{meas1}, for any $a \in A$, since $e(a)=1 \in \RG$, 
$Q_m(a) = (1-q(ma)) \theta_m$ is also in the image of $\psi'$.   Hence so is $\theta_{m,a}$.

Assume now that $A$ is  divisible .  By \propref{meas}, $\RG_{\bdd}[\theta_m: m=1,2,\ldots]$
is contained in the image of $\psi'$.  Moreover if we let 
$$\psi'':  T_A[q^A][\theta, \theta \inv, (1-q(ma_0))^{-1}, (1-(1-\theta \inv)^m) \inv ]_{m=1,2,\ldots}  \to \RG_{l}$$
be the induced homomorphism, then $\psi''$ is surjective.  It follows that  $\psi''$
is an isomorphism.  Let $\psi^*$ be the inverse; restricting back to $\RG$
we obtain the lemma in the  divisible  case.  

In general, define $\psi^*_A$ to be the composition of the natural homomorphism 
$$ \Qq \tensor  \K_{+}^{df}(\volG_A \bdd) \to  \Qq \tensor  \K_{+}^{df}(\volG_{\bA} \bdd) $$
with $\psi^*_{\bA}$.
\eprf

\ssec{The Grothendieck  ring of   $\RV$}

As a step towards the valued field, 
we consider the theory of extensions
$$1 \to \k^* \to \RV \to_{\valr} \G \to 0 $$
 of an ordered  divisible  Abelian group $\G$ (written additively) by the multiplicative group of an 
algebraically closed field.   This is a complete theory; in a saturated model $M$,
the sequence is split, though of course the set of points in a given substructure need not be.  
See \cite{HK} for details.  

 We work over a base structure $A_{\RV}$, which as above is left out of the notation.  
  Let $A$ be the image of $A_{\RV}$ in $\G$.  Let $A_{\RES} = A \meet \RES$
where $\RES = \union_{\g \in \Qq \tensor A} \valr \inv (\g) $.

\

The following specializes Definitions 3.66 and  5.21 of \cite{HK}\,\footnote{The definitions in \cite{HK} are more general in several respects.  In particular
several kinds of resolution on volume forms are considered; here we consider
the type denoted $\vol_\G$ in \cite{HK}.  Since no other volumes are considered, the subscript becomes unnecessary.  Similar results are possible for the other variants.  }.
Define $\Sigma: \G^n \to \G$ by $\Sigma((x_1,\ldots,x_n)) = \sum_{i=1}^n x_i$.  
                                                                                                                                             
\<{defn} \lbl{RVcat} \rm                                                          

1)   $\RV[n]$ is the category of pairs $(U,f)$, with $U$ a 
definable subset  of $\RV^m$ for some $m$, and $f=(f_1,\ldots,f_n): U \to \RV^n$  a finite-to-one map.    A morphism $U \to V$ is a definable bijection $U \to V$.

2)   $\Ob \vol \RV[n] = \Ob \RV[n]$.  A morphism $U \to V$ is a definable
bijection $h: U \to V$ such that for any $u$ we have $\Sigma (f(u))=\Sigma (u)$.

3)    ${\vol \RV \bdd}[m]$ is the full subcategory of ${\vol \RV}[m]$ consisting of objects
whose $\G$-image is contained in 
$[\g,\infty]^m$, for some definable $\g \in \G$.  These will again be referred to as semi-bounded.
   
4)  $\RES[n]$ (respectively $\vol \RES [n]$) is the full subcategory of $RV[n]$ 
(respectively $\vol \RV [n]$) 
whose objects $U$ are contained
in $\RES^m$ for some $m$.  Equivalently, 
such that  $\valr (U)$ is finite.

 \>{defn}

The map $\valr: \RV \to \G$ induces maps $\RV^n \to \G^n$.  If $X, Y$ are
$\G[n]$-isomorphic definable subsets of $\G^n$, then $\valr \inv X, \valr \inv Y$
are definably isomorphic: both $GL_n(\Zz)$ transformations and $A$-translations
obviously lift.  The definition of the category $\G[n]$ was indeed engineered for this.
Hence the pullback $X \mapsto \rv \inv X$ induces a map
\eqn{GtoRV}{ \SG   \G [n] \to \SG \RV[n], \ [X] \mapsto [\valr \inv X]  }

Semi-boundedness is preserved by the pullback ; and also, again by definition, a $\vol \G [n]$-isomorphism lifts to a $\vol \RV[n]$ isomorphism.
Thus we also have
\eqn{vGtoRV}{ \SG   \vol \G \bdd [n] \to \SG \vol \RV \bdd[n], \ [X] \mapsto [\valr \inv X]  }

On the other hand the inclusion induces an obvious map

\eqn{REStoRV}  {\SG \RES [n] \to \SG \RV [n] }
and
\eqn{vREStoRV} {\SG \vol \RES [n] \to \SG \vol \RV [n]}

We obtain  homomorphisms 

\eqn{rv0}{ \SG(\RES[*]) \tensor    \SG({\G}[*]) \to \SG(\RV[*]) }

\eqn{vrv0}{ \SG(\vol \RES [*]) \tensor    \SG({\vol \G \bdd }[*]) \to \SG(\vol \RV \bdd[*]) }

These are shown in \cite{HK1} to be surjective.  
 If $\g \in \G[1]$ is a definable point, then $[\valr \inv (\g)] \in \SG \RES [1]$
has the same image under \eqref{REStoRV} as $\{\g\}_1$ has under \eqref{GtoRV};
and similarly in the measured case.  Thus in both cases the kernel
contains the elements  $ 1 \tensor [\valr \inv (\g)]_1 -  [\g]_1 \tensor 1$,
$\g \in \G$ definable.   By   Corollary 10.3 and Proposition 10.10 of \cite{HK},
{\em these elements generate the kernel } in both cases,  \eqref{REStoRV} and  \eqref{GtoRV}.

\ssec{Bounded definable subsets of $\RV$}
We begin with a description of the Grothendieck ring of two-sided bounded 
definable subsets of $\RV$ in the 
 divisible  case, using \lemref{subring}.      This does not immediately translate to a statement for $\VF$, since 
 the notion of boundedness is not preserved under
arbitrary definable maps.   The results of this subsection will
not be used further on.


\eqref{rv0} induces a homomorphism:

\eqn{rv2}{ \SG(\RES[*])([G_m(k)]_1 \inv ) \tensor   \SG({\G}[*]) ([0]_1 \inv) \to \K(\RV[*]) ( [G_m(\k)]_1 \inv )}
whose kernel is again generated by the elements $ 1 \tensor [\valr \inv (\g)]_1 -  [\g]_1 \tensor 1$,
$\g \in \G$ definable, 
as one can see by multiplying an element of the kernel by a high enough power of $[G_m(\k)]$.

Hence  we have a surjective homomorphism 
\eqn{rv3}{ \SG(\RES[*])([G_m(k)]_1 \inv )_0 \tensor   \SG({\G}[*]) ([0]_1 \inv)_0  \to  \SG(\RV[*]) ( [G_m(\k)]_1 \inv )_0}
whose kernel is generated by the relations $\valr \inv (\g) / [G_m(\k)] =  e(\g)$ (where
$\g \in \Qq \tensor A$, and $e(\g) = [\g]_1  / [0]_1$.)

Note that $\K(\RES[*])([G_m(k)]_1 \inv )_0$ is naturally isomorphic to the direct  limit of the  
$\K(\RES[n])$, where $\K(\RES[n])$ is mapped to $\K(\RES[n+1])$ by the map
$[X] \mapsto [X \times G_m(\k)]$.


\<{defn} $\K^{df}(\RES):= \K(\RES[*])([G_m(k)]_1 \inv )_0$ will be called   the stabilized Grothendieck ring of $\RES$.  Similarly 
$\K^{df}(\Var_F) = \K(\k[*])([G_m(k)]_1 \inv )_0$ and $\K^{df}(\RV) = \K(\RV[*])([G_m(k)]_1 \inv )_0$, and similarly for the semirings.
 \>{defn}

\<{prop} \lbl{I4}  $(\K^{df}(\RES_{A_{\res}}) \tensor \K^{df} (\Gb ) /I \cong \K^{df}(\RV_A \bdd)$  
where
$I$ is the ideal generated by 
$ (\{\frac{\valr \inv (\g)}{[G_m(\k)]} - e(\g): \g \in \Qq \tensor A\}) $
\>{prop}

\prf The homomorphism \eqref{rv0} is compatible with restriction to semi-bounded sets:
  $ \SG(\RES[*]) \tensor    \SG({\G^{\bdd}}[*]) \to \SG(\RV^{\bdd}[*]) $ is surjective 
and has kernel generated by the   elements $1 \tensor [\g] - [\valr \inv (\g) ] \tensor 1$.
Equations \eqref{rv2}, \eqref{rv3}  for semi-bounded sets follow in the same way.
The Proposition follows upon taking additive inverses. 
\eprf

Let $T_A$ denote the symmetric algebra  
$\Qq \oplus (\Qq \tensor A) \oplus Sym^2 (Q \tensor A) \oplus \ldots$.

\<{cor} \lbl{I5}  Assume   $A$ is  divisible , and let $F= A_{\RV} \meet \k$.    Then 
$$ \K^{df}(\RV_A \bdd)   \cong \K^{df}(\Var_{F}) \tensor T_A $$ \>{cor}

\prf  
Assume  $A$ is  divisible .  In this case every definable set $X \subseteq \RES^m$
is definably  isomorphic to a definable subset of a Cartesian power of $\k$, where $\k$ is the
residue field.  So $\K(RES[n])$ reduces to $\K(\k)$, the Grothendieck ring of $F$-varieties.
Moreover  for any definable $\g \in G$, $\valr \inv(\g)$  is definable isomorphic $G_m(\k)$.   Hence in this case the relations in \propref{I4} are redundant, and the tensor product
is valid over $\Qq$.  By \propref{subring}, $ \K_{\Qq}^{df} \Gb \cong T_A$.  The corollary follows.
\eprf

\ssec{The measured Grothendieck ring of $\RV$}

The connection between varieties with forms over the valued field, and the category
$\vol \G[n]$, is mediated by $\vol \RV [n]$.  We now study the dimension-free
Grothendieck ring of this category, incorporating in particular both $\G$ and the residue field.  
 

Let  $F= A_{\RV} \meet \k$ be the base residue field, 
and $\Var_F[n]$ the category of $F$-varieties of dimension $\leq n$.
\eqref{vrv0} can be used to describe $\K^{df}(\vol \RV \bdd )$.  We do this now
in the   case:  $A$ is  divisible .   Recall that  the rings $\K^{df}(\Var_F) , \K^{df}({\volG \bdd}), \K^{df} ({\vol \RV \bdd} )$ are defined with respect to dehomogenizing elements
$[G_m]_1$, $[0]_1$ and $[G_m]_1 \tensor 1 = 1 \tensor [0]_1$ respectively.  

\<{prop}  \lbl{KdfRV} Assume $A$ is  divisible .  Then 

\[    \K^{df} ({\vol \RV \bdd} ) \iso \K^{df}(\Var_F)  \tensor   \K^{df}({\volG \bdd})    \]
 
 \>{prop}
 
\prf  Let    $\K(\Var_F[*]) = \oplus_{n \geq 0} \K(\Var_F[n])$.
  In this case the natural map 
  $$ \K_+(\Var_F [*]) \tensor \K_+(\volG^{fin}[*])   \to \K_+({\vol \RES}[*])$$
  is a surjective homomorphism, with kernel generated by the single relation
  $$R: \ [G_m]_1 \tensor 1 = 1 \tensor [0]_1$$
   
 \eqref{vrv0} simplifies to:

\eq{   \K({\vol \RV \bdd}[*]) \iso \K(\Var_F[*])  \tensor   \K({{\vol \G \bdd}}[*]) / R }

  The proposition follows using  \lemref{gr}.   \eprf

\def\Om{\Omega}

\ssec{The Grothendieck ring of bounded volume forms over valued fields}

  Let $T$ be a $\V$-minimal theory; to simplify notation we will assume $T$ is effective.
See \cite{HK} for the definitions of these notions.  The principal example 
are the theory $ACVF_F$ of algebraically closed valued fields, over a base valued field  $F$ 
with residue field $\bF$ of 
 characteristic $0$.   The reader may take
$T$ to be $ACVF_F$; in this case ``definable'' is the same as ``$F$-semi-algebraic'', and
the category $\Vol_T$ described below is   $\Vol_F$ of the introduction.  
Other examples are analytic expansions of L. Lipshitz and Z. Robinson.

If $V$ is a smooth $n$-dimensional variety, let $\Om V = \bigwedge^{n} T V$, considered
as a variety rather than a vector bundle.
The notion of  a {\em bounded} subset of $V$ and  in the same way as in \cite{serre}, \S 6.1.
 If $X \subseteq V$ is bounded, we consider definable sections $\om: X \to \Om V$
 over $X$; we say $\om$ is bounded if the graph in $\Om V$ is bounded.
 
\<{defn} \lbl{vol} 
 $\Vol_T [n]$ is the category whose objects are pairs $(X,\om)$, with $X$ either empty or 
 a definable bounded Zariski dense subset of  a  smooth $F$-variety $V$
 of dimension $n$,  and 
 $\om: X \to \Om V$   a definable bounded 
section.  
A morphism $ (X,\om) \to (X',\om')$ is a definable bijection $g$ between subsets 
of $X,X'$ whose complement has  dimension $< \dim(V)$, 
such that (away from a set of dimension $< \dim(V)$) $\om = c g ^* \om'$ for some definable function $c$ on $X$ with $\val(c)=0$.
\>{defn}

For $b \in \G$, let $U_b =   \{x:  \val(x)=b\}$.  In particular 
 $U_0 =  \{x:  \val(x)=0\} = \Oo \setminus \Mm$.
  $ \Mm = \{x: \val(x)>0\}$. 
  
$\Vol_T$ is an $\Nn$-graded category, and yields a graded Grothendieck semiring
$\SG(\Vol_T)$.  We take  $e_1=[(U_0,dx)]$, and form the dimension free semiring $\SG^{df}(\Vol_T) = \SG^{df}_{e_1}(\Vol_T) $.  
Let $\K^{df}_{\Qq}(\Vol_T) = \Qq \tensor \SG^{df}(\Vol_T)$.

 To facilitate the comparison to \defref{RVcat}, we  need to compare $\Vol_T$
  to   a  more elementary version.
                                                                                                     
\<{defn} \lbl{VFcat} \rm                                                          

1)   $\VF[n]$ is the category of pairs $(X,f)$, with $X$ a 
definable subset  of $\VF^m $ for some $m$, and $f=(f_1,\ldots,f_n): X \to \VF^n$  a finite-to-one map.    A morphism $X \to Y$ is a definable bijection $X \to Y$.

2)   $\Ob \vol \VF[n] = \Ob \VF[n]$.  A morphism $(X,f) \to (Y,g)$ is a definable
bijection $h: X \to Y$ such that $h^* g^* dx = f^* dx$ away from a variety of dimension $<n$,
where $dx=dx_1 \wedge \ldots \wedge dx_n$ is the standard volume form on $\VF^n$.

3)    ${\vol \VF \bdd}[m]$ is the full subcategory of ${\vol \VF}[m]$ consisting of objects
$(X,f)$ with $f(X)$ bounded.    
 \>{defn}

$\vol \VF$ is dimension-graded, with distinguished element $([U_0]   , Id)$,
and we form $\K_+^{df} \vol \VF$ using the dehomogenizing element $[U_0]$; 
similarly $\K_+^{df} \vol \VF \bdd$ and $\K^{df} \vol \VF \bdd = K(\vol \VF \bdd)^{df}_{[U_0]}$.  

\<{lem} \lbl{vol-df} $\K_+^{df} \vol \VF \cong \K_+^{df} \Vol_T$ canonically;
the isomorphism takes $\K_+^{df} \vol \VF \bdd$ to $\K_+^{df} \vol_T  \bdd$, and induces
an isomorphism   $\K^{df} \vol \VF \cong \K^{df} \Vol_F$.
  \>{lem}

\prf   Let $(X,f) \in \Ob \vol \VF[n]$.  Let $V$ be the Zariski closure
of $X$, and $\om = f^* dx$; this is defined away from a subvariety of $V$ of dimension $< n$.
$(X,f) \mapsto (X,\om)$ is a functor $\VF[n] \to \Vol_T [n]$, inducing an injective
graded semiring homomorphism $\SG \vol \VF  \to \SG \Vol_T$. 
  
 An  element of $\SG \Vol_T [n]$ has the form $[(X,\om)]$ with $X$ a definable
 subset of a smooth affine variety $V \subseteq \VF^{n+l}$, 
 admitting a finite-to-one projection $f: V \to \Aa^n$, and
 $\om (v) = c(v) f^* dx$ for some definable $c: V \to \VF$.  Let
  $Y = \{(x,t) \in V \times \Aa^1: \val(t) = \val(c(x))\}$,  $g(x,t) = (f(x),t)$.
 Then $(X,\om) \times ([U_0],dx) \cong_{\Vol_T} (Y, g^* (dx \wedge dt) ) $
and hence lies in the image of $\vol \VF$. 
 Hence by \lemref{gr3} 
$\SG ^{df} \vol \VF \cong \SG^{df} \Vol_T$ canonically, and so
  $\K^{df} \vol \VF \cong \K^{df} \Vol_F$.
\eprf

  We write 
$\theta_{VF} = 1+ \frac{[\Mm]}{[U_0]}$, and for a definable $b \in \G$ we
write $q_{VF}(b)  = \frac{[U_b]}{[U_0]}$.   These correspond under the canonical
isomorphisms below to the classes $\theta$ and $q(b)$ of $\K^{df}({\vol  \G \bdd})$, and when no confusion can be caused we will omit the subscript.  We assume $\G$ has at least one
definable   element $a_0>0$, and write $\q ^{-m}$ for $q_{VF}(ma_0)$.
 Note that $\q^{-m} = (\q \inv)^m$. 

Write 
$\qq \inv$ for $1- \theta_{VF} \inv \in  \K^{df}( \Vol_T )[\theta_{VF} \inv]$.
So $1 - \qq \inv = \theta_{VF} \inv$

When no confusion can arise, we also write $\q^{-m}$ for  $q(ma_0)$ and
$\qq \inv$ for $1 - \theta \inv$.

Recall  $T_A$ denotes the symmetric algebra  
$\Qq \oplus (\Qq \tensor A) \oplus Sym^2 (Q \tensor A) \oplus \ldots$.

\<{thm} \lbl{vf}  Let $T$ be an effective $\V$-minimal theory.   Let $F$ be the 
field of $\VF$-definable points of $T$, $A = \val(F)$,   $\bA = \Qq \tensor A$, and let
 $0<a_0 \in \bA$.  Then there exists a canonical homomorphism
 $$
  \K^{df}_{\Qq} ( \Vol_T ) \to
  \K^{df}_{\Qq} (\Var_{F^a})[  \qq \inv,   (1-\qq ^{-m}) \inv ]_{m=1,2,\ldots}     ]
    \tensor 
     T_{\bA}[q^{\bA}][  (1-q(ma_0))^{-1}]_{m = 1,2,\ldots}  $$
 
 If  $A$ is  divisible , this induces an isomorphism 
 $$
 \K^{df}_{\Qq} ( \Vol_T )[ \qq \inv, (1-\q ^{-m})^{-1}]_{m}
  \cong  \K^{df}_{\Qq} (\Var_F)[  \qq \inv,   (1-\qq ^{-m}) \inv ]_{m}     ]
    \tensor 
     T_A[q^A][  (1-\q^{-m})^{-1}]_{m}  $$
\>{thm}

\<{remark} (1)  The inverted $1-\qq^{-m}$ on the $\Var_F$ seems to correspond  to 
nothing on the $\Vol_T$-side; see 
 \lemref{theta} and \remref{localization} for an explanation.

(2)  We took  $\K^{df}_{\Qq} (\Var_F)  \{[V]/[G_m^n]: V \in \Var_F, \dim(V) \leq n \}$. 
The localization is by $[G_a]/[G_m]$ and $[G_a^k-[1]_k]/G_m^k, k=1,2,\ldots$.
 
\>{remark}

\<{proof}[Proof of \thmref{vf}]   
Let ${\mathfrak sp}$ be the semiring congruence
on $ \SG {\vol \RV}  $ generated by $([1_\k]_1=[\RVp]_1)$, with the constant $\G$-form $0 \in \G$.
The restriction to $\SG {\vol \RV \bdd}$ is denoted by the same letter, as is the corresponding
ideal of $\K_\Qq {\vol \RV \bdd}$.  
  (The proof of Lemma 8.20 never goes out of the semi-bounded category.) 
 
By \cite{HK} Theorem 8.29 , 
$$\SG( \vol \VF \bdd  [n] ) \cong  \SG({\vol \RV \bdd})/  {\mathfrak sp}$$
  Restricting to $\G$-valued measures as in    (8.5), we obtain an isomorphism
\eq{ \SG( \vol \VF \bdd  [*] )  \cong \SG(\vol  \RV \bdd[*] )  /  {\mathfrak sp} }

If $b = [1_\k]_1 - [\RVp]_1$, 
 this induces a ring isomorphism 
 
\eq{    \K( \vol \VF \bdd [*] )  \cong \K( \vol \RV \bdd  ) [*] / b  }

We take $[G_m(\k)]_1$ as the distinguished element of $\K(\vol \RV  \bdd )[1]$,
and correspondingly the class $[U_0]$ of the annulus $U_0 = \{x:  \val(x)=0$
in $ \K( \Vol_T [1])$; i.e. 
$  \K^{df}  ({\vol \RV \bdd} ) = \K^{df}_{[G_m(\k)]_1} ({\vol \RV \bdd} ) $,
$\K^{df}  ( \Vol_T  ) = \K^{df}_{[U_0]} ( \Vol_T  )$.
 
Let $\xi = \frac{b}{[G_m(\k)]_1}$.   By \lemref{gr2} and \ref{vol-df},
 
 \eqn{vfrv}{ \K^{df}_{\Qq} ( \Vol_T  )= \K^{df}_{\Qq} ( \vol \VF \bdd )  \cong \K^{df}_{\Qq} ( \vol  \RV \bdd_A ) /  \xi}

Thus it suffices to find the canonical homomorphism on $\K^{df}_{\Qq} ( \vol  \RV \bdd_A ) /  \xi $.  This involves work with $\RV$ alone.
At this point we may assume $A$ is  divisible ; the homomorphism in the general case can then be obtained
by composing with the canonical homomorphism $\K^{df}_{\Qq} ( \vol  \RV \bdd_A) \to  \K^{df}_{\Qq} ( \vol  \RV \bdd _{\bf A})$.

By   \propref{KdfRV}, 
\eqn{rvkg}{  \K^{df}_{\Qq} ({\vol \RV \bdd} ) \iso \K^{df}_{\Qq} (\Var_F)  \tensor   \K^{df}_{\Qq}({\vol \G \bdd})    }

Under this isomorphism, $\xi$ corresponds to 

 \eq{\xi_{VF} = \frac{[1]_\k}{[G_m(k)]} \tensor 1 - 1 \tensor \frac{  [\RVp]_1}{q_0}
 = \frac{[1]_\k}{[G_m(k)]} \tensor 1 - 1 \tensor (\theta-1)  }

while $q(ma_0)$ corresponds under the composition of  \eqref{vfrv}, \eqref{rvkg}
to $q_{VF}(ma_0)=\q ^{-m}$, and $\theta$ to $\theta_{VF} $.  `

Hence by \propref{psi}, using $1-\qq \inv = \theta_{VF} \inv$, 

 \eq{ \K^{df}_{\Qq} ( \Vol_T )[ (1-\q ^{-m})^{-1}, \qq \inv]_{m=1,2,\ldots}
 \cong  }
 \eqn{vf-2}{ \K^{df}_{\Qq} (\Var_F)  \tensor 
     T_A[q^A][\theta, \theta \inv, (1-q(ma_0))^{-1}, (1-(1-\theta \inv)^m) \inv ]_{m=1,2,\ldots}  /  \xi_{VF}  }
 
 We can view the relation $\xi_{VF}$ as defining $1 \tensor (\theta-1) = (\qq-1) \inv \tensor 1$
 where $(\qq-1) \inv : = \frac{[1]_k}{G_m(k)]}$.  Then \eqref{vf-2} becomes:
 
  \eq{ \K^{df}_{\Qq} (\Var_F)[ (\qq-1) \inv, \qq \inv,   (1-\qq ^{-m}) \inv ]_{m=1,2,\ldots}     ]
    \tensor 
     T_A[q^A][  (1-q(ma_0))^{-1}]_{ m = 1,2,\ldots} }
As $(\qq -1) \inv = \qq \inv (1-\qq \inv) \inv$,  this term is redundant, so

  \eq{ \K^{df}_{\Qq} ( \Vol_T )[ (1-\q ^{-m})^{-1}, \qq \inv]_{m=1,2,\ldots}
  \cong  }
 \eq{ \K^{df}_{\Qq} (\Var_F)[  \qq \inv,   (1-\qq ^{-m}) \inv ]_{m=1,2,\ldots}     ]
    \tensor 
     T_A[q^A][  (1-q(ma_0))^{-1}]_{m = 1,2,\ldots}}
   \eprf

So far, we always used $[0_\Gamma]_1$ as a dehomogenizing element.   An alternative choice
is $[1_k]_1$; it goes along with  $\Mm$ in $ \K( \Vol_T [1])$ and $[G^{>0}]$ in $\K(vol \G \bdd)$.  
There appears to be a deep duality transposing these choices.  With the latter choice too
one has an analogue of \thmref{qf}, of which we indicate the beginning. 

 Let 
\[ \K^{df'}  ({\vol \RV \bdd} ) = \K ({\vol \RV \bdd}[*] )^{df}_{[1]_1} \]
\[\K^{df'}  ( \Vol_T  ) = \K( \Vol_T  )[*] ^{df}_{[\Mm]}  \]
 \[ \K ^{df'}(\vol \G^{bdd}) = \K(\vol \G^{bdd}) ^{df}_{[RV^{>0}]} \]
 
Also let     $a'= \frac{[G_m]_1 }{[1]_1} \tensor 1 - 1 \tensor \frac{[0]_1}{[RV^{>0}]}$.

\<{lem}  
  \[\K(\vol \RV \bdd_A [*]/ b)^{df}_{[1]_1} \cong  \K^{df'}(\Var_F[*])  \tensor   \K^{df'}({{\vol \G \bdd}}[*]) / a' \]
 \>{lem}
%
 \prf 
   \eqn{vfrv'}{ \K^{df'} ( \Vol_T  )= \K^{df'} ( \vol \VF \bdd )    \cong \K(\vol \RV \bdd_A [*]/ b)^{df}_{[1]_1}     }

 
  
 
In the ring $\K(\vol \RV \bdd_A [*]/ b)$, we have $[1]_1= [RV^{>0}]$, so \lemref{gr} applies.

 As in \propref{KdfRV}, letting $a$ be the ideal generated by $  \ [G_m]_1 \tensor 1 - 1 \tensor [0]_1$, we have  from
 \eqref{vrv0}:

\[  \K({\vol \RV \bdd}[*]) \iso \K(\Var_F[*])  \tensor   \K({{\vol \G \bdd}}[*]) / a  \]
 
 and so
 
\[  \K({\vol \RV \bdd}[*])/b \iso \K(\Var_F[*])  \tensor   \K({{\vol \G \bdd}}[*]) / (a,b)  \]
  
 The statement of the lemma follows from   \lemref{gr}.

 \eprf

If $V$ is a definable subset of a variety  over $F$ and $\om$ a definable volume form, call $(V,\om)$
{\em strictly absolutely integrable} if there exists $(V',\om') \in \Ob \Vol_T$
and a definable bijection $g: V \to V'$ (up to a smaller dimensional set), such that
$\val g^* \om' = \val \om$.  Define $\int_V \om $ to be the image of $[(V',\om')]$
under the homomorphism of \thmref{vf}.  This clearly does not depend on the choice of
$(V',\om')$.

Let $\RR$ be the target ring of \thmref{vf}, and $\int$ the homomorphism.  
$\RR$ admits a natural decreasing $\G$ filtration:  
$$F_\g \RR =   \K^{df}_{\Qq} (\Var_F)[  \qq \inv,   (1-\qq ^{-m}) \inv ]_{m=1,2,\ldots}     ]
    \tensor 
     T_{\bA}[q^{\bA^{>\g}}][  (1-q(ma_0))^{-1}]_{m = 1,2,\ldots}  $$
 where ${\bA^{>\g}} = \{ c \in \bA: c> \g \}$.

 \<{remark}
 Any $(V,\om)$ admits a definable map $c: V \to \G^{\geq 0}$, such that each fiber is 
 strictly absolutely integrable.  Hence so is the inverse image of any bounded subset of $\G$.  Moreover if $V_\g = c \inv (\g)$, then for large $\g$
 $\int_{V_\g} \om = \sum_{i=1}^n r_i P_i(\i(\beta_i \g)) q^{\alpha_i \g}$, with $r_i \in \RR$,
 $P_i \in \Qq[X]$, $\beta_i \in \Qq^m, \alpha_i \in \Qq$.  If all $\alpha_i \geq 0$, and $\alpha_i=0$
 implies $P_i$ is constant, we can call $(V,\om)$  {\em absolutely integrable} and 
 define $\int_V \om = \sum _{\alpha_i =0} r_i P_i$.  This does not depend on the choice of $c$,
 but it is not clear if it is really more general than strict absolute integrability. 
 \>{remark}


 \<{remark}  In \cite{HK} more general volume forms are considered. \rm 
 $   \mu_{\G} \bdd \VF $ is equivalent to the category
   of pairs $(V,\theta)$ with $\theta$ a bounded, bounded support section of the $\G$-bundle  
  $\val_* \bigwedge^{\dim(V)} T V$ induced from the top form bundle via the valuation
  map.   
  If $(V,\om) \in \Vol_F$ then $(V, \val \om)  \in  \mu_{\G} \bdd \VF  $, but the converse need
  not be true.
  
 It is possible to define an integral $\int (V,\theta)$ with values in $\K^{df} \Vol_F$.
 One can easily find definable functions $c: V \to \G$ such that with $V_\g = c \inv (\g)$,
 $(V_\g, \om | V_\g)$ lies in the image of $\Vol_F$.  Then define 
 $\int (V,\theta) = \int _{\g \in \G}  \int_{V_\g} \om | V_\g $.  The expression is well-defined.
 
 However, the dimension-
free  Grothendieck ring   $  \K^{df}_{\Qq} ( \mu_{\G} \bdd \VF )$ is not identical with $\K^{df} (\Vol_F)$.
For instance $\q$ has an  square root in $  \K^{df}_{\Qq} ( \mu_{\G} \bdd \VF )$, namely
$d = [(\{0\},\{\frac{a_0}{2})] / [\{0\},0]$.  We have $d^2 = \q$, 
 as
opposed to the conditional square root $d'=q(\frac{a_0}{2}) \in \K^{df}(\Vol_F)$ which only satisfies 
$(d')^2 = \q e(1/2)$.   Equivalently, the idempotent $e(1/2)$ has a
nontrivial
square root $\frac{d}{d'}$.  

  \>{remark}

 \>{section}

\section{Appendix}

In this appendix we define the Iwahori Hecke algebra of $SL_2 $ over an algebraically closed valued field. We continue to denote by $F$ a valuation field with value group $\Ga$, ring of integers $\cO$ and residue field $\bF$. We denote by $\cO ^{\ga} _o , \cO _{cl} ^{\ga} , \cA ^{\ga}$ the (classes of the) open ball, closed ball and annulus of radius $\ga \in \Ga$. We also denote $q=\cO ^0 _{cl} / \cO ^0 _o$. In particular, $\cA ^0 = (q-1)\cO ^0 _o$. To ease notation, we choose a section $\Ga \to F$ denoted by $\ga \mapsto t^{\ga}$. Note, however, that this is never used in an essential way. 

We denote by $G$ the group $SL_2 (F)$, by $B$ the subgroup of upper triangular matrices, by $N$ the subgroup of unipotent upper triangular matrices and by $A$ the subgroup of diagonal matrices. We will abuse notations and write $G(\cO ), G(\bF )$ etc. for the groups of points of the corresponding algebraic groups. We have a residue map $res:G(\cO )\to G(\bF )$. All integrals over $G$ will be taken with respect to the Haar form on $G$, which is 
\[
dg \bmat  a & b \\\ c & d \emat =\frac{1}{a}da \wedge db \wedge dc
\]
So, for example, the measure of the set of matrices such that $val(a)=\ga _a,val(b)=\ga _b ,val(c)=\ga _c$ is $t^{-\ga _a}\cA ^{\ga _a}\cA ^{\ga _b}\cA ^{\ga _c}=\cA ^0 \cA ^{\ga _b} \cA ^{\ga _c}$.

In order that the convolution makes sense, the field of coefficients will be taken to be a field $E$ together with a ring homomorphism $K^{bdd}(Vol _F)\to E$. By $\propref{psi}$, there is such a field with nontrivial homomorphism.

\begin{defn} A definable function from $G(F)$ to $E$ is a function of the form $f(g)=\sum _{i=0} ^N c_i \phi _i (g)$ where $\phi _i$ are definable functions from $G(F)$ to $K(Var_F)$ and $c_i\in E$. A definable function is called bounded if there is $\ga \in \Ga$ such that $f(g)=0$ unless all entries of $g$ have valuation less than $\ga$.
\end{defn}

\begin{defn} The convolution of two bounded definable functions $f_1 ,f_2$ from $G(F)$ to $K(Var_F)$ is the function $f_1 * f_2 (g)=\int _{h\in G(F)} f_1 (gh^{-1})f_2 (h) dh$, which is easily seen to be a bounded definable function. This definition extends to convolution of bounded definable functions from $G(F)$ to $E$. 
\end{defn}

\begin{rem} We can similarly define bounded definable functions from $\Ga$ to $E$ and convolution of them.
\end{rem}

\begin{defn} The Iwahori subgroup $I\subset G(\cO)$ is the inverse image of $B(\bF )$ under the map $res$. As a vector space, the Iwahori Hecke algebra $\cH$ is the $E$ vector space of bounded definable functions from $G(F)$ to $E$ that are invariant under left and right multiplication by $I$. This is an algebra where the multiplication is convolution of functions.
\end{defn}

A special role will be played by the following $\cH$ module:
\begin{defn} Let $M$ be the right $\cH$ module consisting of bounded definable functions from $G(F)$ to $E$ that are invariant under the left multiplication by $A(\cO )N(F)$ and under the right multiplication by $I$.
\end{defn}

The proof of the following lemmas is standard: 

\begin{lem} Let $g=\bmat x & y \\ z & w \emat$ and $\ga \in \Ga$ be negative. Then
\begin{enumerate}

\item $g\in I\bmat t^{\ga} & 0 \\ 0 & t^{-\ga} \emat I$ iff $val(x)=\ga , val(y)\geq \ga , val(z)>\ga , val(w)\geq \ga$.

\item $g\in I\bmat t^{-\ga} & 0 \\ 0 & t^{\ga} \emat I$ iff $val(x)>\ga , val(y)\geq \ga , val(z)>\ga , val(w)= \ga$.

\item $g\in I\bmat 0 & t^{\ga} \\ t^{-\ga} & 0 \emat I$ iff $val(x)>\ga , val(y)= \ga , val(z)>\ga , val(w)> \ga$.

\item $g\in I\bmat 0 & t^{-\ga} \\ t^{\ga} & 0 \emat I$ iff $val(x)\geq \ga , val(y)\geq \ga , val(z)=\ga , val(w)\geq \ga$.

\end{enumerate}
\end{lem}

\begin{lem} Let $g=\left ( \begin{matrix} x & y \\ z &w \end{matrix} \right )$. Then
\begin{enumerate}

\item If $val(z) \leq val(w)$ then $g\in A(\cO )N \left ( \begin{matrix} 0 & z^{-1} \\ z & 0 \end{matrix} \right ) I$.

\item If $val(z) > val(w)$ then $g\in A(\cO )N \left ( \begin{matrix} w^{-1} & 0 \\ 0 & w \end{matrix} \right ) I$.

\end{enumerate}
\end{lem}

For $\ga \in \Ga$ let $v_{\ga},u_{\ga},S_{\ga},S_{\ga} ^-$ be the characteristic functions of the following double cosets
\[
A(\cO )N \bmat t^{-\ga} & 0 \\ 0 & t^{\ga} \emat I \qq A(\cO )N \bmat 0 & t^{-\ga} \\ t^{\ga} & 0 \emat I \qq I \bmat t^{\ga} & 0 \\ 0 & t^{-\ga} \emat I \qq I \bmat 0 & t^{-\ga} \\ t^{\ga} & 0 \emat I
\]
respectively.

\begin{prop} Let $\ga <0$. Then
\begin{enumerate}

\item $v_0 S_0 = \cA ^0 \cO ^0 _o \cO ^0 _{cl} v_0$

\item $v_0 S_{\ga} = t^{-\ga} \cA ^{\ga} \cO ^{\ga} _o \cO ^{\ga} _{cl}v_{-\ga} + \int \limits _{\ga < \de \leq -\ga} t^{-\ga} \cA ^{\ga} \cA ^{-\de}  \cO ^{\ga} _o u_{\de}$

\item $v_0 S_{-\ga} = \cA ^0 \cO ^0 _o \cO ^0 _{cl}v_{\ga}$

\item $v_0 S^- _0 = \cA ^0 \cO ^0 _o \cO ^0 _{cl}  u_0$

\item $v_0 S^- _{\ga} = t^{-\ga} \cA ^{\ga} \cO ^{0} _o \cO ^{0} _{cl} u_{\ga}$

\item $v_0 S^- _{-\ga} = t^{-\ga} \cA ^{\ga} \cO ^{\ga} _o \cO ^{\ga} _o u_{-\ga} +\int \limits _{\ga < \de < -\ga} t^{-\ga} \cA ^{\ga} \cA ^{-\de} \cO ^{\ga} _ov_{\de}$

\item $u_0 S_0 ^- = \cA ^0 \cO _{cl} ^0 \cO _{cl} ^0 v_0 +\cA ^0 \cA ^0 \cO _{cl} ^0 u_0$

\end{enumerate}
\end{prop}

\begin{proof} We show 2. for example. We first find the coefficients of the $v_{\de}$'s in the convolution. Note that by $I$ invariance, the coefficient of $v_{\de}$ in the convolution equals the value of the convolution at the point $\bmat t^{-\de} & 0 \\ 0 & t^{\de} \emat$. This, in turn, equals to the measure of the set of elements $g\in S_{\ga}$ for which there is $h\in v_0$ such that $gh=\bmat t^{-\de} & 0 \\ 0 & t^{\de} \emat$. Suppose $\bmat x & y \\ z & w \emat \in v_0$ and $\bmat a & b \\ c & d \emat \in S_{\ga}$. If their product is $\bmat t^{-\de} & 0 \\ 0 & t^{\de} \emat$ then
\[
\bmat x & y \\ z & w \emat =\bmat t^{-\de} & 0 \\ 0 & t^{\de} \emat \bmat d & -b \\ -c & a \emat = \bmat t^{-\de}d & -t^{-\de} b \\ -t^{\de}c & t^{\de}a \emat
\]
So $val(t^{\de}a)=0$ and $val(a)=\ga$, hence $\de = -\ga$. The constraints are $val(a)=\ga , val(b),val(d)\geq \ga , val(c)>\ga$, so the coefficient is $t^{-\ga}\cA^{\ga}\cO _o ^{\ga} \cO _{cl} ^{\ga}$. To compute the coefficient of $u_{\de}$, we proceed similarly. Suppose that the product is $\bmat 0 & t^{-\de} \\ t^{\de} & 0 \emat$. Then
\[
\bmat x & y \\ z & w \emat =\bmat 0 & t^{-\de} \\ t^{\de} & 0 \emat \bmat d & -b \\ -c & a \emat = \bmat t^{-\de}c & -t^{-\de} a \\ -t^{\de}d & t^{\de}b \emat
\]
The conditions are $val(a)=\ga , val(b)\geq \ga , val(c)>\ga , val(d) \geq \ga , val(t^{\de}d)>val(t^{\de}b)=0$. This implies $\ga \leq val(b)=-\de$ and $val(d)>val(b)$. We should also have $ad-bc=1$, hence $0=val(ad-bc)\geq \min \{ val(ad),val(bc) \} > \ga -\de$. Hence it is neccessary that $\ga < \de \leq -\ga$. Under this assumption, the conditions are $val(a)=\ga , val(b)= -\de , val(c)>\ga$ (since $val(d)=val(\frac{1}{a} + \frac{bc}{a})\geq \min \{-\ga , -\de \} > \ga$) and the coefficient is $t^{-\ga}\cA ^{\ga} \cA ^{-\de} \cO_o ^{\ga}$.

\end{proof}

We make the following change of base:
\[
e_{\ga}=\frac{1}{\cA ^{-\ga}}v_{-\ga} \QQ f_{\ga}=\frac{1}{\cA ^{\ga}}u_{\ga}
\]
for all $\ga$ and 
\[
S_{\ga}= \cO _{cl} ^0 \cO _{o} ^0 \cA ^{\ga} R_{\ga} \qq S_{-\ga}= \cO _{cl} ^0 \cO _{o} ^0 \cA ^{\ga} R_{-\ga} \qq S_{\ga} ^- = \cO _{cl} ^0 \cO _{o} ^0 \cA ^{\ga} R_{\ga} ^- \qq S_{-\ga} ^-=\cO _{o} ^0 \cO _{o} ^0 \cA ^{\ga} R_{-\ga} ^-
\]
for $\ga<0$. $R_0 ^-$ is defined using the third equality and not the forth. we get 
\begin{cor} \label{cor:mult_tab} Let $\ga <0$. Then
\begin{enumerate}

\item $e_0 R_0 = e_0$

\item $e_0 R_{\ga} = e_{\ga} + \int \limits _{\ga < \de \leq -\ga} \frac{q-1}{q}f_{\de}$

\item $e_0 R_{-\ga} = e_{-\ga}$

\item $e_0 R_{0}^- =f_0$

\item $e_0 R_{\ga}^- = f_{\ga}$

\item $e_0 R_{-\ga} ^- = f_{-\ga} + \int \limits _{\ga < \de <-\ga}(q-1)e_{\de}$

\item $f_0R_0 ^- = qe_0+(q-1)f_0$

\end{enumerate}
\end{cor}

So the transformation $h \mapsto v_0h$ from $\cH$ to $M$ is given by the following block matrix:
\[
\bmat Id & 0 & X & A \\ 0 & Id & 0 & 0 \\ 0 &  0 & Id & 0 \\ B & Y & 0 & Id \emat
\]
where the blocks correspond to the partition $R_{<} , R_{\geq} , R^- _{\leq} , R^- _{>}$ and $e_{<} , e_{\geq} , f_{\leq} , f_{>}$. Here, for example, $A$ is the transformation between two spaces with bases $\{ E_{\ga} \} _{\ga <0}$ and $\{ E'_{\ga} \} _{\ga >0}$ which equals
\[
A E_{\ga}=\frac{q-1}{q}\int _{\de \in (0,-\ga ]} E'_{\de}
\]
This transformation is invertible iff $Id-AB$ is invertible. Now,
\[
(Id-AB)(E_{\ga})=E'_{\ga} -\frac{(q-1)^2}{q}\int \limits_{\eta \in (\ga ,0)}1^{\Ga} _{[\ga ,\eta )}E'_{\eta}
\]
We look for inverse to $Id-AB$ of the form
\[
E'_{\ga} \mapsto E_{\ga} + \int \limits _{\de \in (\ga ,0)} G( \ga - \de ) E_{\de}
\]
The condition on $G$ is that it satisfies
\[
G(z)-\frac{(q-1)^2}{q}1^{\Ga} _{[z,0)} - \frac{(q-1)^2}{q}\int \limits _{w\in (z,0)}G(z-w)1^{\Ga}_{[w,0)} = 0
\]
for every $z<0$. The condition is the same for left and right inverse. There is such a function.
\[
\int \limits _{x\in (\ga ,0)}\cO ^{\ga -x}1_{[x,0)} = \int \limits _{x\in (\ga ,0)} \int \limits _{y\in [x ,0)}\cO ^{\ga -x} = \frac{1}{q-1} \int \limits _{y\in (\ga ,0)} \int \limits _{x\in (\ga ,y]} \cA ^{\ga -x} = \frac{1}{q-1}\int \limits _{y\in (\ga ,0)} \cO _{cl} ^{\ga -y}- \cO _{cl} ^0 = 
\]
\[
\frac{q}{(q-1)^2}\int \limits _{y\in (\ga ,0)}\cA ^{\ga -y} -\frac{q}{q-1}1_{(\ga ,0)} \cO ^0 _{0} = \frac{q}{(q-1)^2}(\cO _o ^{\ga} -\cO _{cl}^0) -\frac{q}{q-1}1_{(\ga ,0)} \cO ^0 _{o} = 
\]
\[
= -\frac{q}{q-1}1_{(\ga ,0)}\cO ^0 _o - \frac{q^2}{(q-1)^2}\cO ^0 _o +\frac{q}{(q-1)^2}\cO ^{\ga} _o
\]
Similarly,
\[
\int \limits _{x\in (\ga ,0)}\cO ^{x-\ga}1_{[x,0)} = \frac{1}{q-1}1_{(\ga ,0)}\cO ^0 _o - \frac{1}{(q-1)^2}\cO ^0 _o +\frac{q}{(q-1)^2}O^{-\ga} _o
\]
From which we see that
\[
G(\ga )=\frac{q}{q^2 -1}\left (\frac{1}{q}\cO ^{\ga} _o -q\cO ^{-\ga} _o \right )
\]
satisfies the equation. 

It follows from the above discussion that $M$ is a rank one free module over $\cH$. In particular, $\cH=End_{\cH}(M)$.

\begin{cor} There is an embedding $\Ga \to \cH$, denoted by $\ga \mapsto \tau _{\ga}$ such that
\[
\tau _{\ga}(v_{\de})=v_{\de- \ga} \QQ \tau _{\ga}(u_{\de})=u_{\de -\ga}
\]
\end{cor}
\begin{proof} $\Ga$ acts on $A(\cO )N \bs G /I=\{ \pm 1\} \ltimes (X_*(A)\otimes \Ga)$ by translations, and hence acts on $M$. The action is $\tau _{\ga}f(g)=f\left ( \bmat t^{-\ga} & 0 \\ 0 & t^{\ga} \emat g \right )$. These transformations are endomorphisms (as left translation commutes with right convolution) so for any $\ga$ there is an unique element $\tau_{\ga}$ acting as the translation. Finally,
\[
(\tau _{\ga}v_{\de})(g)=v_{\de}\left ( \bmat t^{-\ga} & 0 \\ 0 & t^{\ga} \emat g\right ) =1
\]
iff
\[
\bmat t^{-\ga} & 0 \\ 0 & t^{\ga} \emat g\in A(\cO )N\bmat t^{-\de} & 0 \\ 0 & t^{\de} \emat I
\]
iff
\[
g\in \bmat t^{\ga} & 0 \\ 0 & t^{-\ga} \emat A(\cO )N\bmat t^{-\de} & 0 \\ 0 & t^{\de} \emat I = A(\cO )N\bmat t^{\ga} & 0 \\ 0 & t^{-\ga} \emat \bmat t^{-\de} & 0 \\ 0 & t^{\de} \emat = A(\cO )N\bmat t^{\ga -\de} & 0 \\ 0 & t^{\de -\ga} \emat I
\]
iff
\[
v_{\de -\ga}(g)=1
\]

\end{proof}

We have that
\[
\tau _{\ga}(e_{\de})=\frac{\cA ^{-\ga}}{\cA ^0}e_{\de +\ga} \QQ \tau _{\ga}(f_{\de})=\frac{\cA ^{-\ga}}{\cA ^0} f_{\de -\ga}.
\]

This map extends to an embedding of $Fn(\Ga )$ into $\cH$ which is clearly an algebra homomorphism. Denote $T_{\ga}=\frac{\cA ^{\ga}}{\cA ^0} \tau _{\ga}$. The $T_{\ga}$ act as translations on the $e_{\ga},f_{\ga}$'s: $T_{\ga}e_{\de}=e_{\ga+\de},T_{\ga}f_{\de}=f_{\de-\ga}$. We also have $T_{\ga}T_{\de}=T_{\ga+\de}$.

\begin{cor} $(R_0 ^-)^2=(q-1)R_0 ^-+qI$.
\end{cor}

\begin{proof} by computing the action of both sides on $v_0$.
\end{proof}

We let $\hM$ be the set of (definable) functions from $A(\cO )N \bs G /I$ that vanish on $v_{\ga},u_{\ga}$ for $\ga$ negative enough. It is clear that $\hM$ is an $\cH$ module but it is also a $\widehat{\cH}$ module, where $\widehat{\cH}$ is the obvious completion of $\cH$. We define $\sI:M \to \hM$ by 
\[
\sI \varphi (g)=\int _N \varphi (wng)dn
\]
Where $N=\left \{ \bmat 1 & x \\ 0 & 1 \emat \right \}$ is the unipotent upper triangular matrices  and $w=\bmat 0 & 1 \\ -1 & 0 \emat$ is the nontrivial element of the Weyl group. 

\begin{lem} $\sI$ is a well defined homomorphism of $\cH$ modules. We have
\begin{enumerate}

\item $\sI T_{\ga}=T_{-\ga} \sI$

\item $\sI e_0=\cO ^0 _{o} f_0 + \cA ^0 \int _{(0,\infty)} T_{\ga} e_0$.

\item $\sI f_0 =\cO _{cl} ^0 e_0 +\cA ^0 \int _{[0,\infty )}T_{\ga} f_0$

\item $\sI (e_0+f_0)=(\cO _{cl} ^0 T_0 + \cA ^0 \int _{\ga \in (0,\infty )} T_{\ga} )(e_0+f_0)$

\end{enumerate}
\end{lem}

\begin{proof} (1) Denote $W_{\ga}=\bmat t^{-\ga} & 0 \\ 0 & t^{\ga} \emat$. Then $wW_{\ga}=W_{-\ga}w$ and $W_{\ga}\bmat 1 & x \\ 0 & 1 \emat W_{-\ga}=\bmat 1 & t^{-2\ga} x \\ 0 & 1 \emat$.
\[
\sI \tau _{\ga} f(g)=\int _{N} (\tau _{\ga} f)(wng)dn =\int _{N}f(W_{\ga} wng)dn = \int _{N}f(wW_{-\ga}nW_{\ga}W_{-\ga}g)dn =
\]
The change of coordinates $m= W_{-\ga}nW_{\ga}$ satisfies $dm=\frac{\cA ^{2\ga}}{\cA ^0}dn$
\[
=\frac{\cA ^{-2\ga}}{\cA ^0}\int _{N}f(wmW_{-\ga}g)dm = \frac{\cA ^{-2\ga}}{\cA ^0}\tau _{-\ga} \sI f (g)
\]
So
\[
\sI T_{\ga}=\frac{\cA ^{\ga}}{\cA ^0}\sI \tau _{\ga}=\frac{\cA ^{\ga}\cA ^{-2\ga}}{\cA ^0 \cA ^0}\tau _{-\ga}\sI =T_{-\ga}\sI
\]
(2) Let $n=\bmat 1 & x \\ 0 & 1 \emat$. Then $wn=\bmat 0 & 1 \\ -1 & -x \emat$. To compute the coefficient of $v_{\ga}$, suppose
\[
wn\bmat t^{-\ga} & 0 \\ 0 & t^{\ga} \emat = \bmat 0 & t^{\ga} \\ -t^{-\ga} & -xt^{\ga} \emat \in ANI
\]
Then $-\ga = val(t^{-\ga}) >val(xt^{\ga})=0$. Hence $\ga <0$ and the measure of $x$'s that contribute is $\cA ^{-\ga}$. If, on the other hand, 
\[
wn\bmat 0 & t^{-\ga} \\ -t^{\ga} & 0 \emat = \bmat t^{\ga} & 0 \\ -xt^{\ga} & -t^{-\ga} \emat \in ANI
\]
Then $val(xt^{\ga})>val(t^{-\ga})=0$, hence $\ga =0$ and the measure of $x$'s that contributes is $\cO _{o} ^0$. Hence $\sI (v_0)=\cO ^0 _o u_0 + \int _{(-\infty ,0)}\cA ^{-\ga} v_{\ga}$, so $\sI e_0=\cO ^0 _{o} f_0 + \cA ^0 \int _{(0,\infty)} T_{\ga} e_0$.

(3) Similarly, assume
\[
wn\bmat t^{-\ga} & 0 \\ 0 & t^{\ga} \emat = \bmat 0 & t^{\ga} \\ -t^{-\ga} & -xt^{\ga} \emat \in ANwI
\]
Then $0=val(t^{-\ga}) \leq val(xt^{\ga})$, so $\ga =0$ and the measure of $x$'s is $\cO _{cl} ^0$. If, on the other hand,
\[
wn\bmat 0 & t^{-\ga} \\ -t^{\ga} & 0 \emat = \bmat t^{\ga} & 0 \\ -xt^{\ga} & -t^{-\ga} \emat \in ANwI
\]
Then $0=val(xt^{\ga})\leq val(t^{-\ga})=-\ga$, so $\ga \leq 0$ and the measure of $x$'s is $\cA ^{-\ga}$. Hence $\sI u_0=\cO _{cl} ^0 v_0 +\int _{(-\infty ,0]}\cA ^{-\ga}u_{\ga}$. This implies that $\sI f_0 =\cO _{cl} ^0 e_0 +\cA ^0 \int _{[0,\infty )}T_{\ga} f_0$.
(4) follows from (2) and (3).
\end{proof}

Note that $\sI$ does not preserve $M$. However, we claim that the operator $J_b=(1-T_{b})\sI$ preserves $M$ for every $b \in \Ga$. Take for example $b>0$. By computing the action on $e_0$ we see that
\[
J_b=\cO _o (1 -T_b)R_0 ^- +\cA ^0 \int _{(0,b]} T_{\ga}.
\]

Fix $a\in \Ga$ and let $b>0$ be smaller in absolute value. Using $J_b T_a=T_{-a}J_b$ and the last equality we get
\[
(1-T_b)\cO ^0 _{o} R_0 ^- T_a + \cA ^0 \int _{(a,a+b]}T_{\ga} = (1-T_b)\cO _{o} ^0 T_{-a}R_0 ^- + \cA ^0 \int _{(-a,-a+b]} T_\ga
\]
and so if $a>0$,
\[
(1-T_b)\cO _o(R_0 ^- T_a^- -T_{-a}R_0 ^-)=\cA ^0 (\int _{(-a,-a+b]}T_{\ga} -\int _{(a,a+b]}T_{\ga})=(1-T_b)\cA ^0 \int _{(-a,a]}T_{\ga}
\]
and if $a<0$,
\[
(1-T_b)\cO _o(R_0 ^- T_a^- -T_{-a}R_0 ^-)=-(1-T_b)\cA ^0 \int _{(a,-a]}T_{\ga}
\]

\begin{lem} The element $1-T_b$ does not annihilate non zero elements of $\cH$.
\end{lem}

\begin{proof} Suppose $X\in \cH$ is non zero. We can view $X$ as a definable function from $\{ \pm 1\}\ltimes \Ga$ to $E$. The support of $X$ is a definable set, hence there is a supremum $\ga$ for it. Let $\ep \in \Ga$ be positive and smaller than $b$ such that $X(\ga -\ep )\neq 0$. Then $(1-T_b)X(\ga+b-\ep )\neq 0$, so $(1-T_b )X\neq 0$.
\end{proof}

\begin{cor} (Bernstein's presentation) Every element in $\cH$ is of the form $\int _{\Ga} f_1(\ga )T_{\ga} + \int _{\Ga} f_{w}(\ga )T_{\ga}R_0 ^-$. Multiplication is defined by being $\Ga$-additive and the relations
\[
R_0 ^-T_{a}=T_{-a}R_0 ^- +(q-1)\int _{(-a,a]}T_{\ga} 
\]
and
\[
(R_0 ^- -q)(R_0 ^- +1)=0
\]
\end{cor}

\begin{prop} The center of $\cH$ consists of all elements of the form $\int _{\Ga} f(\ga)(T_{\ga}+T_{-\ga})$.
\end{prop}

\begin{proof} Denote by $L$ the algebra (or space) generated by the $T_{\ga}+T_{-\ga}$. Clearly, $L$ is contained in the center. On the other hand, every element in $\cH$ can be uniquely written as a combination of elements of the form $T_{\ga}+T_{-\ga},T_{\ga}-T_{-\ga}, (T_{\ga}+T_{-\ga})R_0 ^-,(T_{\ga}-T_{\ga})R_0 ^-$ (note that $T_{\ga}e_0 =e_{\ga}$ and $R_0 ^- e_0=f_0$). Every one of those subspaces is $L$ invariant and they are linearly independent.
\end{proof}

\begin{cor} The algebra $\cH$ is finite over its center.
\end{cor}

\<{thebibliography}{serre}

\bibitem{HK}  E. Hrushovski, D. Kazhdan, Integration over valued fields
 
 \bibitem{HKP} Haines, Thomas J., Kottwitz, Robert E., Prasad, Amritanshu, Iwahori Hecke algebras, math.RT/0309168
 
  \bibitem{fesenko}   
   Invitation to higher local fields. Papers from the conference held in MŸnster, August 29--September 5, 1999. Edited by Ivan Fesenko and Masato Kurihara. Geometry and Topology Monographs, 3. Geometry and Topology Publications, Coventry, 2000. front matter+304 pp. (electronic)
   
   \bibitem{kageyama-fujita}  M. Kageyama, M. Fujita, Grothendieck rings of o-minimal expansions
of ordered Abelian groups,   arXiv:math.LO/0505331 v1 16 May 2005 

\bibitem{Lee} Lee, Kyu-Hwan, Iwahori-Hecke algebras of $SL_2$ over 2-dimensional local fields. Preprint 2006
    
  \bibitem{marikova}   Marikova, J., 
  MA thesis,  Charles University, Prague 2003; 
  Geometric properties of semilinear and semibounded sets, preprint.

     \bibitem{serre}   Serre, Jean-Pierre Lectures on the Mordell-Weil theorem. Translated from the French and edited by Martin Brown from notes by Michel Waldschmidt. Aspects of Mathematics, E15.   Vieweg  Braunschweig, 1989. 
\>{thebibliography}
\end{document}